\PassOptionsToPackage{linktocpage}{hyperref}
\documentclass[review]{elsarticle}
\usepackage[fleqn]{amsmath}
\usepackage{amssymb,latexsym,amsmath}
\usepackage{overpic}
\usepackage{color}
    \usepackage[hyperindex,breaklinks]{hyperref}
\usepackage{graphicx}
\usepackage{amsfonts}
\usepackage{epstopdf}
\usepackage{epsfig}
\graphicspath{{../Plots/}}

\newtheorem{thm}{Theorem}[section]
\newtheorem{defn}{Definition}[section]

\newtheorem{lem}{Lemma}[section]

\newtheorem{exmp}{Example}[section]
\numberwithin{equation}{section}

\begin{document}
\begin{frontmatter}
\title{A high-order nodal discontinuous Galerkin method for  nonlinear fractional
Schr\"{o}dinger type equations}
\author[]{Tarek Aboelenen}
\ead{tarek.aboelenen@aun.edu.eg}




\address{Department of Mathematics, Assiut University, Assiut 71516, Egypt}
\begin{abstract}
 We propose  a nodal discontinuous Galerkin  method  for  solving the nonlinear
 Riesz space fractional Schr\"{o}dinger equation and  the strongly coupled nonlinear Riesz space fractional Schr\"{o}dinger equations. These problems have been expressed
as a system of low order differential/integral equations. Moreover, we prove, for both problems,  $L^{2}$ stability  and optimal order of convergence $O(h^{N+1})$, where $h$ is space step size and $N$ is polynomial degree. Finally, the performed numerical experiments confirm the optimal order of convergence.




\end{abstract}
\end{frontmatter}
\section{Introduction}
In this paper we develop a nodal discontinuous Galerkin method to solve the generalized nonlinear fractional Schr\"{o}dinger  equation
\begin{equation}\label{sch1vn}
\begin{split}
&i\frac{\partial u}{\partial t}- \lambda_{1}(-\Delta)^{\frac{\alpha}{2}}u+ \lambda_{2}f(|u|^{2})u=0,\\
&u(x,0) = u_{0}(x),
\end{split}
\end{equation}
and the strongly coupled nonlinear fractional Schr\"{o}dinger equations
\begin{equation}\label{sch2}
\begin{split}
&i\frac{\partial u}{\partial t}- \lambda_{1}(-\Delta)^{\frac{\alpha}{2}}u+ \varpi_{1}u+\varpi_{2} v+ \lambda_{2}f(|u|^{2},|v|^{2})u=0,\\
&i\frac{\partial v}{\partial t}- \lambda_{3}(-\Delta)^{\frac{\alpha}{2}}v+\varpi_{2}u+ \varpi_{1}v+\lambda_{4} g(|u|^{2},|v|^{2})v=0,\\
&u(x,0) = u_{0}(x),\\
&v(x,0) =  v_{0}(x),
\end{split}
\end{equation}
and homogeneous boundary conditions. $f(u)$ and $g(u)$ are arbitrary (smooth) nonlinear real functions and  $\lambda_{i}$, $i=1,2,3,4$  are a real constants,
$\varpi_{1}$ is normalized birefringence constant and $\varpi_{2}$ is the linear coupling parameter which accounts for the effects that arise from twisting and elliptic deformation of the fiber \cite{cai2010multisymplectic}. Notice that the assumption of homogeneous boundary conditions is for simplicity only and is not essential: the method
can be easily designed for nonhomogeneous boundary conditions. The  fractional Laplacian $-(-\Delta)^{\frac{\alpha}{2}}$, which can be defined using Fourier analysis as  \cite{muslih2010riesz, yang2010numerical}

$$-(-\Delta)^{\frac{\alpha}{2}}u(x,t)=\mathcal{F}^{-1}(|\xi|^{\alpha}\hat{u}(\xi,t))$$
where $\mathcal{F}$ is the Fourier transform. Equation \eqref{sch1vn} can be viewed as a generalization of the classical nonlinear  Schr\"{o}dinger  equation. During the last decade, it has arisen as a suitable model in many application areas, such as  fluid dynamics, nonlinear optics, and plasma physics \cite{bialynicki1979gaussons,bullough1979solitons,cowan1986quasi}. It was first introduced by Laskin \cite{laskin2000fractional,laskin2000fractionalcf},  who derived fractional Schr\"{o}dinger equation with Riesz space-fractional derivative includes a space fractional derivative of order $\alpha\,\,(1 < \alpha < 2)$ instead of the Laplacian in the classical Schr\"{o}dinger equation, and obtained its by replacing Brownian trajectories in Feynman path integrals (corresponding to the classical Schr\"{o}dinger equation) by the L\'{e}vy flights.
  It is generally difficult to give the explicit forms of the analytical solutions of nonlinear fractional Schr\"{o}dinger equation, thus the construction of numerical methods becomes very important.
In recent years, developing various numerical algorithms for solving nonlinear fractional Schr\"{o}dinger  equation has received much attention. For the
time-fractional Schr\"{o}dinger equation, Wei et al.\cite{wei2012analysis} presented and analyzed an implicit fully discrete local discontinuous Galerkin (LDG) finite element method for solving the time-fractional Schr\"{o}dinger equation.  Hicdurmaza and  Ashyralyev   presented  stability analysis for a first order difference scheme
applied to a nonhomogeneous  multidimensional time fractional Schr\"{o}dinger differential equation. For the space-fractional Schr\"{o}dinger equation, Wang and Huang \cite{wang2015energy} studied  an energy conservative Crank-Nicolson difference scheme for nonlinear Riesz space-fractional Schr\"{o}dinger equation. Yang \cite{yang2016class} proposed a class of linearized energy-conserved finite difference schemes for nonlinear space-fractional Schr\"{o}dinger equation. Galerkin finite element method for nonlinear fractional Schr\"{o}dinger equations were considered \cite{li2016galerkin}.  Amore et.al. \cite{amore2010collocation} developed
the collocation method for fractional quantum mechanics.\\
The strongly coupled nonlinear Schr\"{o}dinger system \eqref{sch2} arise in many physical fields, especially in
in fluid mechanics, solid state physics and plasma waves and for two interacting nonlinear packets in a dispersive and conservative system, see, e.g.,\cite{benney1967propagation, yang1997classification,ran2016conservative} and reference therein.
 When $\alpha=2$,  it represents the integer-order strongly coupled equations, and a number of conservative schemes for such case have been proposed \cite{sonnier2005strong,wang2008numerical,ismail2007linearly}. When $\varpi_{1}=\varpi_{2}=0$, this system becomes the weakly coupled   nonlinear fractional Schr\"{o}dinger equations considered in \cite{wang2013crank,li2016galerkin}  and reference therein. Ran and Zhang \cite{ran2016conservative} proposed a conservative difference scheme for solving the strongly
coupled nonlinear fractional Schr\"{o}dinger equations.
A numerical study based on an implicit fully discrete LDG for the time-fractional coupled Schr\"{o}dinger systems is presented \cite{wei2012numerical}. To the best of our knowledge, however,
the LDG method,
which is an important approach to solve partial differential equations and fractional partial differential equations, has not been
considered for  the nonlinear Schr\"{o}dinger equation and the coupled nonlinear  Schr\"{o}dinger equations  with the Riesz space fractional derivative. Compared with finite difference methods, it has the
advantage of greatly facilitates the handling of complicated geometries and elements
of various shapes and types, as well as the treatment of boundary conditions.\\
The   LDG method is a well-established method for classical
conservation laws \cite{Hesthaven:2007:NDG:1557392,Bernardo,ref1}.  For application of the method to fractional problems, Mustapha and McLean \cite{Mustapha2011,mustapha2013superconvergence}  have developed and analyzed discontinuous Galerkin methods for time fractional diffusion and wave equations.  Xu and Hesthaven \cite{doi:10.1137/130918174}  proposed a LDG method for  fractional convection-diffusion
equations. They proved stability and optimal order of convergence $N+1$
for the fractional diffusion problem when polynomials of degree $N$, and an order of convergence of $N+\frac{1}{2}$ is established for the general fractional convection-diffusion problem with general monotone flux for the nonlinear term.  Aboelenen and El-Hawary \cite{cann} proposed a high-order nodal discontinuous Galerkin  method for  a linearized fractional  Cahn-Hilliard equation. They proved stability and optimal order of convergence $N+1$
for the linearized fractional Cahn-Hilliard  problem. Here we propose LDG method for problems \eqref{sch1vn}-\eqref{sch2}  with the Riesz space fractional derivative of order $\alpha$ $(1 < \alpha< 2)$.
For $1 < \alpha< 2$, it is conceptually
similar to a fractional derivative with an order between $1$ and $2$. We rewrite the fractional operator as a composite of first order derivatives and a fractional integral  and convert the nonlinear fractional Schr\"{o}dinger equation and the strongly coupled nonlinear fractional  Schr\"{o}dinger  equations  into a system of low order equations. This allows us to apply the LDG method.\\
The outline of this paper is as follows. In section \ref{sc1}, we introduce some basic definitions and recall a few central results. In section \ref{sc2}, we derive the discontinuous
Galerkin formulation for the  nonlinear fractional Schr\"{o}dinger equation. In section \ref{sc3}, we prove a theoretical result of $L^{2}$ stability for the nonlinear case as well as an error estimate for the linear case. In section \ref{sc40} we present a local discontinuous Galerkin method for the strongly coupled nonlinear  fractional Schr\"{o}dinger equations and give a theoretical result of $L^{2}$ stability  for the nonlinear case and an error estimate for the linear case in section \ref{sc4}. Section \ref{sc5} presents some numerical examples to  illustrate the efficiency of the scheme. A few concluding remarks are offered in section \ref{sc6}.

  \section{Preliminary definitions}\label{sc1}
 We introduce some preliminary definitions of fractional calculus, see, e.g.,\cite{miller1993introduction} and associated functional setting for the subsequent numerical
schemes and theoretical analysis.
 \subsection{Liouville-Caputo Fractional Calculus}
 The left-sided and right-sided Riemann-Liouville integrals of order $\alpha$, when $0 < \alpha < 1$, are defined, respectively, as
 \begin{equation}\label{1}
\big({}^{\,\,RL}_{-\infty}\mathcal{I}_{x}^{\alpha}f\big)(x)=\frac{1}{\Gamma(\alpha)}\int_{-\infty}^{x}
\frac{f(s)ds}{(x-s)^{1-\alpha}}, \quad x > -\infty,
\end{equation}
and
\begin{equation}\label{1111}
\big({}^{RL}_{\,\,\,\,x}\mathcal{I}_{\infty}^{\alpha}f\big)(x)=
\frac{1}{\Gamma(\alpha)}\int_{x}^{\infty}\frac{f(s)ds}{(s-x)^{1-\alpha}}, \quad x < \infty,
\end{equation}
where $\Gamma$ represents the Euler Gamma function. The corresponding inverse operators, i.e., the left-sided and
right-sided  fractional derivatives of order $\alpha$, are then defined based on \eqref{1} and \eqref{1111}, as
\begin{equation}\label{2}
\big({}^{\,\,RL}_{-\infty}\mathcal{D}_{x}^{\alpha}f\big)(x)=\frac{d}{dx}\big({}^{\,\,RL}_{-\infty}
\mathcal{I}_{x}^{1-\alpha}f\big)(x)=\frac{1}{\Gamma(1-\alpha)}
\frac{d}{dx}\int_{-\infty}^{x}\frac{f(s)ds}{(x-s)^{\alpha}}, \quad x > -\infty,
\end{equation}
and
\begin{equation}\label{3}
\big({}^{RL}_{\,\,\,\,x}\mathcal{D}_{\infty}^{\alpha}f\big)(x)
=\frac{-d}{dx}\big(^{RL}_{\,\,\,\,x}\mathcal{I}_{\infty}
^{1-\alpha}f\big)(x)=\frac{1}{\Gamma(1-\alpha)}
\bigg(\frac{-d}{dx}\bigg)\int_{x}^{\infty}\frac{f(s)ds}{(s-x)^{\alpha}}, \quad x < \infty.
\end{equation}
This allows for the definition of the left and right Riemann-Liouville fractional derivatives of order $\alpha$ $ (n-1 <\alpha<n),\,\,n\in \mathbb{N}$ as
\begin{equation}\label{2}
\big({}^{\,\,RL}_{-\infty}\mathcal{D}_{x}^{\alpha}f\big)(x)=\bigg(\frac{d}{dx}\bigg)^{n}\big({}^{\,\,RL}_{-\infty}
\mathcal{I}_{x}^{n-\alpha}f\big)(x)=\frac{1}{\Gamma(n-\alpha)}
\bigg(\frac{d}{dx}\bigg)^{n}\int_{-\infty}^{x}\frac{f(s)ds}{(x-s)^{-n+1+\alpha}}, \quad x > -\infty,
\end{equation}
and
\begin{equation}\label{3}
\big({}^{RL}_{\,\,\,\,x}\mathcal{D}_{\infty}^{\alpha}f\big)(x)=\bigg(\frac{-d}{dx}\bigg)^{n}\big(^{RL}_{\,\,\,\,x}
\mathcal{I}_{\infty}
^{n-\alpha}f\big)(x)=\frac{1}{\Gamma(n-\alpha)}
\bigg(\frac{-d}{dx}\bigg)^{n}\int_{x}^{\infty}\frac{f(s)ds}{(s-x)^{-n+1+\alpha}}, \quad x < \infty.
\end{equation}
Furthermore, the corresponding left-sided and right-sided  Caputo derivatives of order $\alpha$ $ (n-1 <\alpha<n)$ are obtained as
\begin{equation}\label{4n}
\big({}^{\,\,\,\,\,\,C}_{-\infty}\mathcal{D}_{x}^{\alpha}f\big)(x)=\bigg({}^{RL}_{-\infty}\mathcal{I}_{x}^{n-\alpha}
\frac{d^{n}f}{dx^{n}}\bigg)(x)
=\frac{1}
{\Gamma(n-\alpha)}\int_{-\infty}^{x}\frac{f^{(n)}(s)ds}{(x-s)^{-n+1+\alpha}}, \quad x > -\infty,
\end{equation}
and
\begin{equation}\label{5b}
\big({}^{C}_{\,x}\mathcal{D}_{\infty}^{\alpha}f\big)(x)=(-1)^{n}\bigg({}^{RL}_{\,\,\,\,x}
\mathcal{I}_{\infty}^{n-\alpha}\frac{d^{n}f}{dx^{n}}
\bigg)(x)=
\frac{1}{\Gamma(n-\alpha)}\int_{x}^{\infty}\frac{(-1)^{n}f^{(n)}(s)ds}{(s-x)^{-n+1+\alpha}}, \quad x < \infty.
\end{equation}
The Riesz fractional derivative is defined as
\begin{equation}\label{n5}
\frac{\partial^{\alpha}}{\partial |x|^{\alpha}}u(x,t)=-(-\Delta)^{\frac{\alpha}{2}}u(x,t)=-\frac{
{}^{\,\,\,\,\,\,C}_{-\infty}\mathcal{D}_{x}^{\alpha}u(x,t)+{}^{C}_{\,x}\mathcal{D}_{\infty}^{\alpha}u(x,t)}
{2\cos\big(\frac{\pi\alpha}{2}\big)}.
\end{equation}
If $\alpha<0$, the fractional Laplacian becomes the fractional integral operator. In this case, for any $0 <\mu<1$, we define
\begin{equation}\label{nc5}
\Delta_{-\mu/2}u(x)=-\frac{
{}^{\,\,\,\,\,\,C}_{-\infty}\mathcal{D}_{x}^{-\mu}u(x)
+{}^{C}_{\,x}\mathcal{D}_{\infty}^{-\mu}u(x)}{2\cos\big(\frac{\pi(2-\mu)}{2}\big)}=\frac{
{}^{\,\,\,\,\,\,C}_{-\infty}\mathcal{D}_{x}^{-\mu}u(x)
+{}^{C}_{\,x}\mathcal{D}_{\infty}^{-\mu}u(x)}{2\cos\big(\frac{\pi\mu}{2}\big)}=\frac{
{}^{RL}_{-\infty}\mathcal{I}_{x}^{-\mu}u(x)
+{}^{RL}_{\,\,\,\,x}
\mathcal{I}_{\infty}^{-\mu}u(x)}{2\cos\big(\frac{\pi\mu}{2}\big)}.
\end{equation}
When $1 <\alpha<2$, using \eqref{4n}, \eqref{5b} and \eqref{nc5}, we can rewrite the fractional Laplacian in
the following form:
\begin{equation}\label{n5}
-(-\Delta)^{\frac{\alpha}{2}}u(x)=\Delta_{\frac{(\alpha-2)}{2}}\bigg(\frac{d^{2}u(x)}{dx^{2}}\bigg).
\end{equation}
To carry out the analysis, we introduce the appropriate fractional spaces.
 \begin{defn}(left fractional space \cite{Ervin_variationalformulation}). We define the seminorm
\begin{equation}\label{7}
|u|_{J_{L}^{\alpha}(\mathbb{R})}=\big\|{}^{RL}_{\,x_{L}}\mathcal{D}_{x}^{\alpha}u\big\|_{L^{2}(\mathbb{R})}.
\end{equation}
and the norm
\begin{equation}\label{7}
\|u\|_{J_{L}^{\alpha}(\mathbb{R})}=(|u|_{J_{L}^{\alpha}(\mathbb{R})}^{2}+\|u\|_{L^{2}(\mathbb{R})}^{2})^{\frac{1}{2}},
\end{equation}
and let $J_{L}^{\alpha}(\mathbb{R})$ denote the closure of $C_{0}^{\infty}(\mathbb{R})$ with respect to $\|.\|_{J_{L}^{\alpha}(\mathbb{R})}$.\\
\end{defn}
\begin{defn} (right fractional space \cite{Ervin_variationalformulation}). We define the seminorm
\begin{equation}\label{7}
|u|_{J_{R}^{\alpha}(\mathbb{R})}=\big\|{}^{RL}_{\,\,\,\,x}
\mathcal{D}_{x_{R}}^{\alpha}u\big\|_{L^{2}(\mathbb{R})},
\end{equation}
and the norm
\begin{equation}\label{7}
\|u\|_{J_{R}^{\alpha}(\mathbb{R})}=(|u|_{J_{R}^{\alpha}(\mathbb{R})}^{2}+\|u\|_{L^{2}(\mathbb{R})}^{2})^{\frac{1}{2}},
\end{equation}
and let $J_{R}^{\alpha}(\mathbb{R})$ denote the closure of $C_{0}^{\infty}(\mathbb{R})$ with respect to $\|.\|_{J_{R}^{\alpha}(\mathbb{R})}$.
\end{defn}
\begin{defn} (symmetric fractional space \cite{Ervin_variationalformulation}). We define the seminorm
\begin{equation}\label{7}
\|u\|_{J_{S}^{\alpha}(\mathbb{R})}=\big|\big({}^{RL}_{\,x_{L}}\mathcal{D}_{x}^{\alpha}u,{}^{RL}_{\,\,\,\,x}
\mathcal{D}_{x_{R}}^{\alpha}u\big)_{L^{2}(\mathbb{R})}\big|^{\frac{1}{2}},
\end{equation}
and the norm
\begin{equation}\label{7}
\|u\|_{J_{S}^{\alpha}(\mathbb{R})}=\big(|u|_{J_{S}^{\mu}(\mathbb{R})}^{2}+\|u\|_{L^{2}(\mathbb{R})}^{2}\big)^{\frac{1}{2}}.
\end{equation}
and let $J_{S}^{\alpha}(\mathbb{R})$ denote the closure of $C_{0}^{\infty}(\mathbb{R})$ with respect to $\|.\|_{J_{S}^{\alpha}(\mathbb{R})}$.
\end{defn}
\begin{lem}\label{lk}(see \cite{Ervin_variationalformulation}). For any $0 <s<1$, the fractional integral satisfies the following
property:
\begin{equation}\label{7}
({}^{\,\,RL}_{-\infty}\mathcal{I}_{x}^{s}u,{}^{RL}_{\,\,\,\,x}\mathcal{I}_{\infty}^{s}u)_{\mathbb{R}}=
\cos(s\pi)|u|_{J_{L}^{-s}(\mathbb{R})}^{2}=
\cos(s\pi)|u|_{J_{R}^{-s}(\mathbb{R})}^{2}.
\end{equation}
\end{lem}
\begin{lem}\label{lkb} For any $0 <\mu<1$, the fractional integral satisfies the following
property:
\begin{equation}\label{7}
(\Delta_{-\mu}u,u)_{\mathbb{R}}=
|u|_{J_{L}^{-\mu}(\mathbb{R})}^{2}=
|u|_{J_{R}^{-\mu}(\mathbb{R})}^{2}.
\end{equation}
\end{lem}
Generally, we consider the problem in a bounded domain instead of $\mathbb{R}$. Hence,
we restrict the definition to the domain $\Omega = [a, b]$.
\begin{defn} Define the spaces $J_{R,0}^{\alpha}(\Omega),J_{L,0}^{\alpha}(\Omega),J_{S,0}^{\alpha}(\Omega)$ as the closures of
$C_{0}^{\infty}(\Omega)$ under their respective norms.
\end{defn}
\begin{lem}\label{lg}
(fractional Poincar$\acute{e}$-Friedrichs, \cite{Ervin_variationalformulation}). For $u \in J_{L,0}^{\alpha}(\Omega)$ and $\alpha \in \mathbb{R}$, we have
\begin{equation}\label{7}
\|u\|_{L^{2}(\Omega)}\leq C |u|_{J_{L,0}^{\alpha}(\Omega)},
\end{equation}
and for $u \in J_{R,0}^{\alpha}(\Omega)$, we have
\begin{equation}\label{7}
\|u\|_{L^{2}(\Omega)}\leq C |u|_{J_{R,0}^{\alpha}(\Omega)}.
\end{equation}
\end{lem}
\begin{lem}\label{lga234} (See \cite{Kilbas:2006:TAF:1137742})
For any $0 <\mu<1$, the fractional integration operator ${}^{RL}_{\,\,\,\,a}\mathcal{I}_{x}^{\mu}$ is bounded in $L^{2}(\Omega)$:
\begin{equation}\label{7}
\|{}^{RL}_{\,\,\,\,a}\mathcal{I}_{x}^{\mu}u\|_{L^{2}(\Omega)}\leq K \|u\|_{L^{2}(\Omega)}.
\end{equation}
The fractional integration operator ${}^{RL}_{\,\,\,\,x}\mathcal{I}_{b}^{\mu}$ is bounded in $L^{2}(\Omega)$:
\begin{equation}\label{7}
\|{}^{RL}_{\,\,\,\,x}\mathcal{I}_{b}^{\mu}u\|_{L^{2}(\Omega)}\leq K \|u\|_{L^{2}(\Omega)}.
\end{equation}
\end{lem}

\begin{lem}\label{lga2}
The fractional integration operator $\Delta_{-\mu}$ is bounded in $L^{2}(\Omega)$:
\begin{equation}\label{7}
\|\Delta_{-\mu}u\|_{L^{2}(\Omega)}\leq K \|u\|_{L^{2}(\Omega)}.
\end{equation}
\end{lem}
\textbf{Proof}. Combining Lemma \ref{lga234} with \eqref{nc5}, we obtain the result.

\section{ LDG method for  nonlinear fractional
Schr\"{o}dinger equation}\label{sc2}
Let us consider nonlinear fractional Schr\"{o}dinger equation. To obtain a
high order discontinuous Galerkin scheme for the fractional derivative, we rewrite the
fractional derivative as a composite of first order derivatives and a fractional integral
to recover the equation to a low order system. However, for the first order system,
alternating fluxes are used. We introduce three variables $e,r,s$ and set
\begin{equation}\label{1a}
\begin{split}
&e=\Delta_{(\alpha-2)/2}r, \quad r=\frac{\partial}{\partial x}s,\quad s=\frac{\partial}{\partial x}u,
\end{split}
\end{equation}
then, the nonlinear fractional
Schr\"{o}dinger  problem can be rewritten as
\begin{equation}\label{1bchcf}
\begin{split}
&i\frac{\partial u}{\partial t}+\lambda_{1}e+ \lambda_{2}f(|u|^{2})u=0,\\
&e=\Delta_{(\alpha-2)/2}r,\quad r=\frac{\partial}{\partial x}s,\quad s=\frac{\partial}{\partial x}u.\\
\end{split}
\end{equation}
For actual numerical implementation, it might be more efficient if we decompose the complex function
$u(x,t)$ into its real and imaginary parts by writing
\begin{equation}\label{1b}
\begin{split}
u(x,t)=p(x,t)+iq(x,t),
\end{split}
\end{equation}
where $p$, $q$ are real functions. Under the new notation, the problem \eqref{1bchcf} can be written as
\begin{equation}\label{1bx}
\begin{split}
&\frac{\partial p}{\partial t}+\lambda_{1}e+ \lambda_{2}f(p^{2}+q^{2})q=0,\\
&e=\Delta_{(\alpha-2)/2}r,\quad r=\frac{\partial}{\partial x}s,\quad s=\frac{\partial}{\partial x}q,\\
&\frac{\partial q}{\partial t}-\lambda_{1}l- \lambda_{2}f(p^{2}+q^{2})p=0,\\
&l=\Delta_{(\alpha-2)/2}w,\quad w=\frac{\partial}{\partial x}z,\quad z=\frac{\partial}{\partial x}p.\\
\end{split}
\end{equation}

We consider problems posed on the physical domain $\Omega$ with boundary $\partial\Omega$ and assume that this domain is well approximated by the computational domain $\Omega_{h}$. We consider a nonoverlapping element $D^{k}$ such that
\begin{equation}\label{82aa1}
\Omega \simeq \Omega_{h}=\bigcup_{k=1}^{K}D^{k}.
\end{equation}
Now we introduce the broken Sobolev space for any real number $r$
\begin{equation}\label{82aa1}
H^{r}(\Omega_{h})=\{v\in L^{2}(\Omega):\forall k=1,2,....K,v|_{D^{k}}\in H^{r}(D^{k})\}.
\end{equation}
We define the local inner product and $L^{2}(D^{k})$ norm
\begin{equation}\label{82aa1}
(u,v)_{D^{k}}=\int_{D^{k}}uvdx,\quad \|u\|^{2}_{D^{k}}=(u,u)_{D^{k}},
\end{equation}
as well as the global broken inner product and norm
\begin{equation}\label{82aa1}
(u,v)_{\Omega_{h}}=\sum_{k=1}^{K}(u,v)_{D^{k}},\quad \|u\|^{2}_{L^{2}(\Omega_{h})}=\sum_{k=1}^{K}(u,u)_{D^{k}}.
\end{equation}
To complete the LDG scheme, we introduce  the numerical flux.\\
The numerical traces $( p, q,s,z)$ are defined on interelement faces as the
alternating fluxes \cite{doi:10.1137/S0036142997316712,ref1}
\begin{equation}\label{fl}
\begin{split}
p^{*}_{k+\frac{1}{2}}=p^{-}_{k+\frac{1}{2}},\,\,s^{*}_{k+\frac{1}{2}}
=s^{+}_{k+\frac{1}{2}},\,\,\,\,
q^{*}_{k+\frac{1}{2}}=q^{-}_{k+\frac{1}{2}},
z^{*}_{k+\frac{1}{2}}=z^{+}_{k+\frac{1}{2}}.
\end{split}
\end{equation}
Note that we can also choose
\begin{equation}\label{1a}
\begin{split}
p^{*}_{k+\frac{1}{2}}=p^{+}_{k+\frac{1}{2}},\,\,s^{*}_{k+\frac{1}{2}}
=s^{-}_{k+\frac{1}{2}},\,\,\,\,
q^{*}_{k+\frac{1}{2}}=q^{+}_{k+\frac{1}{2}},
z^{*}_{k+\frac{1}{2}}=z^{-}_{k+\frac{1}{2}}.
\end{split}
\end{equation}
For  simplicity we  discretize  the  computational domain $\Omega$ into $K$ non-overlapping elements,  $D^{k}=[x_{k-\frac{1}{2}},x_{k+\frac{1}{2}}]$, $\Delta x_{k}=x_{k+\frac{1}{2}}-x_{k-\frac{1}{2}}$ and $k = 1,...,K$. Let $p_{h}, q_{h}, e_{h},l_{h},r_{h},s_{h},w_{h},z_{h}\in V_{k}^{N}$ be the approximation of $p,q,e,l,r,s,w,z$ respectively, where  the
approximation space  is defined as
\begin{equation}\label{82aa1}
V_{k}^{N}=\{v:v_{k}\in\mathbb{P}(D^{k}), \, \forall D^{k}\in \Omega\},
\end{equation}
where $\mathbb{P}(D^{k})$ denotes the set of polynomials of degree up to $N$  defined  on  the  element  $D^{k}$.
 We define  local discontinuous Galerkin scheme as follows: find $p_{h}, q_{h}, e_{h},l_{h},r_{h},s_{h},w_{h},z_{h}\in V_{k}^{N}$, such that for all test functions $\vartheta_{1},\beta_{1},\phi,\varphi,\chi,\beta_{2},\psi,\zeta\in V_{k}^{N}$,
\begin{equation}\label{1bch1nb}
\begin{split}
&\big(\frac{\partial p_{h}}{\partial t},\vartheta_{1}\big)_{D^{k}}+\lambda_{1}\big(e_{h},\vartheta_{1}\big)_{D^{k}}+ \lambda_{2}\big(f(p_{h}^{2}+q_{h}^{2})q_{h},\vartheta_{1}\big)_{D^{k}}=0,\\
&\big(e_{h},\beta_{1}\big)_{D^{k}}=\big(\Delta_{(\alpha-2)/2}r_{h},\beta_{1}\big)_{D^{k}},\\
&\big(r_{h},\phi\big)_{D^{k}}=\big(\frac{\partial}{\partial x}s_{h},\phi\big)_{D^{k}},\\
&\big(s_{h},\varphi\big)_{D^{k}}=\big(\frac{\partial}{\partial x}q_{h},\varphi\big)_{D^{k}},\\
&\big(\frac{\partial q_{h}}{\partial t},\chi\big)_{D^{k}}-\lambda_{1}\big(l_{h},\chi\big)_{D^{k}}- \lambda_{2}\big(f(p_{h}^{2}+q_{h}^{2})p_{h},\chi\big)_{D^{k}}=0,\\
&\big(l_{h},\beta_{2}\big)_{D^{k}}=\big(\Delta_{(\alpha-2)/2}w_{h},\beta_{2}\big)_{D^{k}},\\
&\big(w_{h},\psi\big)_{D^{k}}=\big(\frac{\partial}{\partial x}z_{h},\psi\big)_{D^{k}},\\
&\big(z_{h},\zeta\big)_{D^{k}}=\big(\frac{\partial}{\partial x}p_{h},\zeta\big)_{D^{k}}.\\
\end{split}
\end{equation}
Applying integration by parts to \eqref{1bch1nb}, and replacing the fluxes at the interfaces
by the corresponding numerical fluxes, we obtain
\begin{equation}\label{1bch2}
\begin{split}
&\big( (p_{h})_{t},\vartheta_{1}\big)_{D^{k}}+\lambda_{1}\big(e_{h},\vartheta_{1}\big)_{D^{k}}+ \lambda_{2}\big(f(p_{h}^{2}+q_{h}^{2})q_{h},\vartheta_{1}\big)_{D^{k}}=0,\\
&\big(e_{h},\beta_{1}\big)_{D^{k}}=\big(\Delta_{(\alpha-2)/2}r_{h},\beta_{1}\big)_{D^{k}},\\
&\big(r_{h},\phi\big)_{D^{k}}=-\big(s_{h},\phi_{x}\big)_{D^{k}}+\big(n.s^{*}_{h},\phi\big)_{\partial D^{k}},\\
&\big(s_{h},\varphi\big)_{D^{k}}=-\big(q_{h},\varphi_{x}\big)_{D^{k}}+\big(n.q^{*}_{h},\varphi\big)_{\partial D^{k}},\\
&\big((q_{h})_{t},\chi\big)_{D^{k}}-\lambda_{1}\big(l_{h},\chi\big)_{D^{k}}- \lambda_{2}\big(f(p_{h}^{2}+q_{h}^{2})p_{h},\chi\big)_{D^{k}}=0,\\
&\big(l_{h},\beta_{2}\big)_{D^{k}}=\big(\Delta_{(\alpha-2)/2}w_{h},\beta_{2}\big)_{D^{k}},\\
&\big(w_{h},\psi\big)_{D^{k}}=-\big(z_{h},\psi_{x}\big)_{D^{k}}+\big(n.z^{*}_{h},\psi\big)_{\partial D^{k}},\\
&\big(z_{h},\zeta\big)_{D^{k}}=-\big(p_{h},\zeta_{x}\big)_{D^{k}}+\big(n.p^{*}_{h},\zeta\big)_{\partial D^{k}}.\\
\end{split}
\end{equation}
\section{ Stability and error estimates}\label{sc3}
 In the following we discuss stability and accuracy of the proposed scheme, for the nonlinear fractional Schr\"{o}dinger   problem.
\subsection{ Stability analysis }  In order to carry out the analysis of the LDG scheme, we have the following results.
\begin{thm}\label{tt4h} ($L^{2}$ stability).
The semidiscrete scheme \eqref{1bch2} is   stable, and $
\|u_{h}(x,T)\|_{\Omega_{h}}\leq c\|u_{0}(x)\|_{\Omega_{h}}$ for any $T>0$.
\end{thm}
\textbf{Proof.}
Set $(\vartheta_{1},\beta_{1},\phi,\varphi,\chi,\beta_{2},\psi,\zeta)=(p_{h},-r_{h}+e_{h},p_{h},-z_{h},q_{h},l_{h}-w_{h},-q_{h},s_{h})$ in \eqref{1bch2}, and consider the integration by parts formula $\big(u,\frac{\partial r}{\partial x}\big)_{D^{k}}+\big(r,\frac{\partial u}{\partial x}\big)_{D^{k}}=[ur]_{_{x_{k-\frac{1}{2}}}}^{x_{k+\frac{1}{2}}}$, we get
\begin{equation}\label{1bch1}
\begin{split}
&\big((p_{h})_{t},p_{h}\big)_{D^{k}}+\big((q_{h})_{t},q_{h}\big)_{D^{k}}
+\big(e_{h},e_{h}\big)_{D^{k}}+\big(l_{h},l_{h}\big)_{D^{k}}
+\big(\Delta_{(\alpha-2)/2}w_{h},w_{h}\big)_{D^{k}}+\big(\Delta_{(\alpha-2)/2}r_{h},r_{h}\big)_{D^{k}}\\
&=\big(\Delta_{(\alpha-2)/2}w_{h},l_{h}\big)_{D^{k}}
+\big(\Delta_{(\alpha-2)/2}r_{h},e_{h}\big)_{D^{k}}-\big( r_{h},p_{h}\big)_{D^{k}}+\big( w_{h},q_{h}\big)_{D^{k}}+\big( e_{h},r_{h}\big)_{D^{k}}+\big( l_{h},w_{h}\big)_{D^{k}}\\
&\quad-\lambda_{1}\big( e_{h},p_{h}\big)_{D^{k}}+\lambda_{1}\big( l_{h},q_{h}\big)_{D^{k}}+\theta(s_{h},p_{h})-\theta(q_{h},z_{h}),\\
\end{split}
\end{equation}
with entropy fluxes
\begin{equation}\label{b4}
\begin{split}
\theta(u,v)=\big(n.u^{*},v\big)_{\partial D^{k}}+\big(n.v^{*},u\big)_{\partial D^{k}}-\big(n.u,v\big)_{\partial D^{k}}.\\
\end{split}
\end{equation}
Employing Young's inequality and  Lemma \ref{lga2}, we obtain
\begin{equation}\label{1bch1}
\begin{split}
&\big((p_{h})_{t},p_{h}\big)_{D^{k}}+\big((q_{h})_{t},q_{h}\big)_{D^{k}}
+\big(e_{h},e_{h}\big)_{D^{k}}+\big(l_{h},l_{h}\big)_{D^{k}}
+\big(\Delta_{(\alpha-2)/2}w_{h},w_{h}\big)_{D^{k}}+\big(\Delta_{(\alpha-2)/2}r_{h},r_{h}\big)_{D^{k}}\\
&\leq c_{4}\|p_{h}\|^{2}_{L^{2}(D^{k})}+c_{3}\|q_{h}\|^{2}_{L^{2}(D^{k})}
+c_{6}\|w_{h}\|^{2}_{L^{2}(D^{k})}+c_{5}\|r_{h}\|^{2}_{L^{2}(D^{k})}+c_{1}\|e_{h}\|^{2}_{L^{2}(D^{k})}
\\
&+c_{2}\|l_{h}\|^{2}_{L^{2}(D^{k})}
+\theta(s_{h},p_{h})-\theta(q_{h},z_{h}).\\
\end{split}
\end{equation}
Recalling Lemma \ref{lg}, provided $c_{i},\,\,i=1,2,3,4$ are sufficiently small such that $c_{i}\leq1$, we obtain that
\begin{equation}\label{1bch1n}
\begin{split}
&\big((p_{h})_{t},p_{h}\big)_{D^{k}}+\big((q_{h})_{t},q_{h}\big)_{D^{k}}
\leq \|p_{h}\|^{2}_{L^{2}(D^{k})}+\|q_{h}\|^{2}_{L^{2}(D^{k})}
+\theta(s_{h},p_{h})-\theta(q_{h},z_{h}),\\
\end{split}
\end{equation}
we notice that, with the definition \eqref{fl} of the numerical fluxes
and with simple algebraic manipulations and summing over all elements \eqref{1bch1n}, we easily obtain
\begin{equation}\label{cdgh}
\begin{split}
\sum_{k=1}^{K}(\theta(s_{h},p_{h})-\theta(q_{h},z_{h}))=0. \\
\end{split}
\end{equation}
This implies that
\begin{equation}\label{1bch1}
\begin{split}
\big((p_{h})_{t},p_{h}\big)_{L^{2}(\Omega_{h})}+\big((q_{h})_{t},q_{h}\big)_{L^{2}(\Omega_{h})}
\leq
\|p_{h}\|^{2}_{\Omega_{h}}+\|q_{h}\|^{2}_{\Omega_{h}}.\\
\end{split}
\end{equation}
Hence
\begin{equation}\label{cd}
\begin{split}
\frac{1}{2}\frac{d}{dt}\|u_{h}(x,t)\|^{2}_{\Omega_{h}}\leq \|u(x,t)\|^{2}_{\Omega_{h}}.\\
\end{split}
\end{equation}
Employing Gronwall's inequality, we obtain $
\|u_{h}(x,T)\|_{\Omega_{h}}\leq c\|u_{0}(x)\|_{\Omega_{h}}$.$\quad\Box$

\subsection{Error estimates}
We consider the linear fractional Schr\"{o}dinger equation
\begin{equation}\label{1bchj}
\begin{split}
&i\frac{\partial u}{\partial t}-\lambda_{1}(-\Delta) ^{\frac{\alpha}{2}}u+\lambda_{2}u=0.\\
\end{split}
\end{equation}
It is easy to verify that the exact solution of the above  \eqref{1bchj} satisfies
\begin{equation}\label{1bch2j}
\begin{split}
&\big( p_{t},\vartheta_{1}\big)_{D^{k}}+\lambda_{1}\big(e,\vartheta_{1}\big)_{D^{k}}+ \lambda_{2}\big(q,\vartheta_{1}\big)_{D^{k}}=0,\\
&\big(e,\beta_{1}\big)_{D^{k}}=\big(\Delta_{(\alpha-2)/2}r,\beta_{1}\big)_{D^{k}},\\
&\big(r,\phi\big)_{D^{k}}=-\big(s,\phi_{x}\big)_{D^{k}}+\big(n.s^{*},\phi\big)_{\partial D^{k}},\\
&\big(s,\varphi\big)_{D^{k}}=-\big(q,\varphi_{x}\big)_{D^{k}}+\big(n.q^{*},\varphi\big)_{\partial D^{k}},\\
&\big(q_{t},\chi\big)_{D^{k}}-\lambda_{1}\big(l,\chi\big)_{D^{k}}-\lambda_{2} \big(p,\chi\big)_{D^{k}}=0,\\
&\big(l,\beta_{2}\big)_{D^{k}}=\big(\Delta_{(\alpha-2)/2}w,\beta_{2}\big)_{D^{k}},\\
&\big(w,\psi\big)_{D^{k}}=-\big(z,\psi_{x}\big)_{D^{k}}+\big(n.z^{*},\psi\big)_{\partial D^{k}},\\
&\big(z,\zeta\big)_{D^{k}}=-\big(p,\zeta_{x}\big)_{D^{k}}+\big(n.p^{*},\zeta\big)_{\partial D^{k}}.\\
\end{split}
\end{equation}
Subtracting  \eqref{1bch2j}, from  the linear fractional Schr\"{o}dinger equation \eqref{1bch2}, we have the following error equation
\begin{equation}\label{1bch2k2}
\begin{split}
\big( &(p-p_{h})_{t},\vartheta_{1}\big)_{D^{k}}+\big((q-q_{h})_{t},\chi\big)_{D^{k}}
-\big(\Delta_{(\alpha-2)/2}(r-r_{h}),\beta_{1}\big)_{D^{k}}
-\big(\Delta_{(\alpha-2)/2}(w-w_{h}),\beta_{2}\big)_{D^{k}}\\
&+\big(s-s_{h},\phi_{x}\big)_{D^{k}} +\big(q-q_{h},\varphi_{x}\big)_{D^{k}}+\big(z-z_{h},\psi_{x}\big)_{D^{k}}
+\big(p-p_{h},\zeta_{x}\big)_{D^{k}}\\
&+ \lambda_{2}\big(q-q_{h},\vartheta_{1}\big)_{D^{k}}- \lambda_{2}\big(p-p_{h},\chi\big)_{D^{k}}
+\big(r-r_{h},\phi\big)_{D^{k}}+\big(s-s_{h},\varphi\big)_{D^{k}}
+\big(l-l_{h},\beta_{2}\big)_{D^{k}}\\
&+\big(e-e_{h},\beta_{1}\big)_{D^{k}}+\big(w-w_{h},\psi\big)_{D^{k}} +\big(z-z_{h},\zeta\big)_{D^{k}}-\big(n.(s-s_{h})^{*},\phi\big)_{\partial D^{k}}-\lambda_{1}\big(l-l_{h},\chi\big)_{D^{k}}\\
&+\lambda_{1}\big(e-e_{h},\vartheta_{1}\big)_{D^{k}}-\big(n.(q-q_{h})^{*},\varphi\big)_{\partial D^{k}}
-\big(n.(z-z_{h})^{*},\psi\big)_{\partial D^{k}}-\big(n.(p-p_{h})^{*},\zeta\big)_{\partial D^{k}}=0.\\
\end{split}
\end{equation}

For the error estimate, we define special projections, $\mathcal{P}^{-}$ and $\mathcal{P}^{+}$ into $V_{h}^{k}$.  For
all the elements, $D^{k}$, $k = 1, 2, ... , K$ are defined to satisfy
\begin{equation}\label{prh}
\begin{split}
&(\mathcal{P}^{+}u-u,v)_{D^{k}}=0,\quad\forall  v\in\mathbb{P}_{N}^{k}(D^{k}),\quad \mathcal{P}^{+}u(x_{k-\frac{1}{2}})=u(x_{k-\frac{1}{2}}),\\
&(\mathcal{P}^{-}u-u,v)_{D^{k}}=0,\quad\forall  v\in\mathbb{P}_{N}^{k-1}(D^{k}),\quad \mathcal{P}^{-}u(x_{k+\frac{1}{2}})=u(x_{k+\frac{1}{2}}).\\
\end{split}
\end{equation}
Denoting
\begin{equation}\label{91h}
\begin{split}
&\pi=\mathcal{P}^{-}p-p_{h},\quad \pi^{e}=\mathcal{P}^{-}p-p,\quad \epsilon=\mathcal{P}^{+}r-r_{h},\quad \epsilon^{e}=\mathcal{P}^{+}r-r,\quad\phi_{1}=\mathcal{P}^{+}e-e_{h},\quad\phi_{1}^{e}=\mathcal{P}^{+}e-e,\\
& \tau=\mathcal{P}^{+}s-s_{h},\quad\tau^{e}=\mathcal{P}^{+}s-s,\quad \sigma=\mathcal{P}^{-}q-q_{h}, \quad \sigma^{e}=\mathcal{P}^{-}q-q,\quad\phi_{2}=\mathcal{P}^{+}l-l_{h},\quad\phi_{2}^{e}=\mathcal{P}^{+}l-l,  \\
&\varphi=\mathcal{P}^{+}w-w_{h},\quad \varphi^{e}=\mathcal{P}^{+}w-w,\quad \vartheta=\mathcal{P}^{+}z-z_{h},\quad \vartheta^{e}=\mathcal{P}^{+}z-z.
\end{split}
\end{equation}
For the special projections mentioned above, we have, by the standard approximation
theory \cite{Ciarlet:2002:FEM:581834}, that
\begin{equation}\label{sth}
\begin{split}
&\|\mathcal{P}^{+}u(.)-u(.)\|^{2}_{L^{2}(\Omega_{h})}\leq  Ch^{N+1},\\
&\|\mathcal{P}^{-}u(.)-u(.)\|^{2}_{L^{2}(\Omega_{h})}\leq  Ch^{N+1},\\
\end{split}
\end{equation}
where here and below $C$ is a positive constant (which may have a different value in each
occurrence) depending solely on u and its derivatives but not of  $h$.\\
 \begin{lem}\label{lm1h}
\begin{equation}\label{91ch}
\begin{split}
&\big(\frac{\partial \pi}{\partial t},\pi\big)_{\Omega_{h}}+\big(\frac{\partial \sigma}{\partial t},\sigma\big)_{\Omega_{h}}+\big(\Delta_{(\alpha-2)/2}\epsilon,\epsilon\big)_{\Omega_{h}}
+\big(\Delta_{(\alpha-2)/2}\varphi,\varphi\big)_{\Omega_{h}}+\big(\phi_{1},\phi_{1}\big)_{\Omega_{h}}
+\big(\phi_{2},\phi_{2}\big)_{\Omega_{h}}\\
&=Q_{1}+Q_{2}+Q_{3}+Q_{4}, \\
\end{split}
\end{equation}
where
\begin{subequations}\label{91chjm}
\begin{align}
&Q_{1}=-\big(\epsilon,\pi\big)_{\Omega_{h}}
+\big(\varphi,\sigma\big)_{\Omega_{h}}+\big(\Delta_{(\alpha-2)/2}\epsilon,\phi_{1}\big)_{\Omega_{h}}
+\big(\Delta_{(\alpha-2)/2}\varphi,\phi_{2}\big)_{\Omega_{h}}\\
&\quad\quad\quad-\lambda_{1}\big(\phi_{1},\pi\big)_{\Omega_{h}}
+\lambda_{1}\big(\phi_{2},\sigma\big)_{\Omega_{h}} +\big(\phi_{2},\varphi\big)_{\Omega_{h}}
+\big(\phi_{1},\epsilon\big)_{\Omega_{h}}\\
&Q_{2}=\big(\tau^{e},\pi_{x}\big)_{\Omega_{h}} -\big(\sigma^{e},\vartheta_{x}\big)_{\Omega_{h}}-\big(\vartheta^{e},\sigma_{x}\big)_{\Omega_{h}}
+\big(\pi_{h},\tau_{x}\big)_{\Omega_{h}}
+\big(\vartheta^{e},\tau\big)_{\Omega_{h}}-\big(\tau^{e},\vartheta\big)_{\Omega_{h}} ,\\
&Q_{3}=\big( (\pi^{e})_{t},\pi\big)_{\Omega_{h}}+\big((\sigma^{e})_{t},\sigma\big)_{\Omega_{h}}
+\big(\phi_{2}^{e},\phi_{2}-\varphi\big)_{\Omega_{h}}+\big(\phi_{1}^{e},\phi_{1}-\epsilon\big)_{\Omega_{h}}+\lambda_{2}\big(\sigma^{e},\pi\big)_{\Omega_{h}}\\
&\quad\quad- \lambda_{2}\big(\pi^{e},\sigma\big)_{\Omega_{h}}+\big(\epsilon^{e},\pi\big)_{\Omega_{h}}
-\big(\varphi^{e},\sigma\big)_{\Omega_{h}}-\big(\Delta_{(\alpha-2)/2}\epsilon^{e},\phi_{1}-\epsilon\big)_{\Omega_{h}}
-\big(\Delta_{(\alpha-2)/2}\varphi^{e},\phi_{2}-\varphi\big)_{\Omega_{h}}\\
&\quad\quad+\lambda_{1}\big(\phi_{1}^{e},\pi\big)_{\Omega_{h}}
-\lambda_{1}\big(\phi_{2}^{e},\sigma\big)_{\Omega_{h}}
,\\
&Q_{4}=-\sum_{k=1}^{K}((\tau^{e})^{+}[\pi])_{k+\frac{1}{2}}
+\sum_{k=1}^{K}((\sigma^{e})^{-}[\vartheta])_{k+\frac{1}{2}}+
\sum_{k=1}^{K}((\vartheta^{e})^{+}[\sigma])_{k+\frac{1}{2}}-\sum_{k=1}^{K}((\pi^{e})^{-}[\tau])_{k+\frac{1}{2}}. \\
\end{align}
\end{subequations}

\end{lem}
\textbf{Proof.} From the Galerkin orthogonality \eqref{1bch2k2}, we get
\begin{equation}\label{1bch2k}
\begin{split}
\big( &(\pi-\pi^{e})_{t},\vartheta_{1}\big)_{D^{k}}+\big((\sigma-\sigma^{e})_{t},\chi\big)_{D^{k}}
-\big(\Delta_{(\alpha-2)/2}(\epsilon-\epsilon^{e}),\beta_{1}\big)_{D^{k}}
-\big(\Delta_{(\alpha-2)/2}(\varphi-\varphi^{e}),\beta_{2}\big)_{D^{k}}\\
&+\big(\tau-\tau^{e},\phi_{x}\big)_{D^{k}} +\big(\sigma-\sigma^{e},\varphi_{x}\big)_{D^{k}}+\big(\vartheta-\vartheta^{e},\psi_{x}\big)_{D^{k}}
+\big(\pi-\pi_{h},\zeta_{x}\big)_{D^{k}}\\
&+\lambda_{2} \big(\sigma-\sigma^{e},\vartheta_{1}\big)_{D^{k}}- \lambda_{2}\big(\pi-\pi^{e},\chi\big)_{D^{k}}
+\big(\epsilon-\epsilon^{e},\phi\big)_{D^{k}}+\big(\tau-\tau^{e},\varphi\big)_{D^{k}}
+\big(\phi_{2}-\phi_{2}^{e},\beta_{2}\big)_{D^{k}}\\
&+\big(\phi_{1}-\phi_{1}^{e},\beta_{1}\big)_{D^{k}}+\big(\varphi-\varphi^{e},\psi\big)_{D^{k}} +\big(\vartheta-\vartheta^{e},\zeta\big)_{D^{k}}+\lambda_{1}\big(\phi_{1}-\phi_{1}^{e},\vartheta_{1}\big)_{D^{k}}-\lambda_{1}\big(\phi_{2}-\phi_{2}^{e},\chi\big)_{D^{k}}
\\
&-\big(n.(\tau-\tau^{e})^{*},\phi\big)_{\partial D^{k}}-\big(n.(\sigma-\sigma^{e})^{*},\varphi\big)_{\partial D^{k}}
-\big(n.(\vartheta-\vartheta^{e})^{*},\psi\big)_{\partial D^{k}}-\big(n.(\pi-\pi^{e})^{*},\zeta\big)_{\partial D^{k}}=0.\\
\end{split}
\end{equation}

We take the test functions
\begin{equation}\label{91h}
\begin{split}
\vartheta_{1}=\pi,\quad\beta_{1}=\phi_{1}-\epsilon,\quad \phi=\pi,\quad \varphi=-\vartheta,\quad \chi=\sigma,\quad\beta_{2}=\phi_{2}-\varphi,\quad \psi=-\sigma,\quad \zeta=\tau,
\end{split}
\end{equation}
 we obtain
\begin{equation}\label{1bch2k}
\begin{split}
\big( &(\pi-\pi^{e})_{t},\pi\big)_{D^{k}}+\big((\sigma-\sigma^{e})_{t},\sigma\big)_{D^{k}}
-\big(\Delta_{(\alpha-2)/2}(\epsilon-\epsilon^{e}),\phi_{1}-\epsilon\big)_{D^{k}}
-\big(\Delta_{(\alpha-2)/2}(\varphi-\varphi^{e}),\phi_{2}-\varphi\big)_{D^{k}}\\
&+\big(\tau-\tau^{e},\pi_{x}\big)_{D^{k}} -\big(\sigma-\sigma^{e},\vartheta_{x}\big)_{D^{k}}-\big(\vartheta-\vartheta^{e},\sigma_{x}\big)_{D^{k}}
+\big(\pi-\pi_{h},\tau_{x}\big)_{D^{k}}\\
&+\lambda_{2} \big(\sigma-\sigma^{e},\pi\big)_{D^{k}}-\lambda_{2} \big(\pi-\pi^{e},\sigma\big)_{D^{k}}
+\big(\epsilon-\epsilon^{e},\pi\big)_{D^{k}}-\big(\tau-\tau^{e},\vartheta\big)_{D^{k}}
+\big(\phi_{2}-\phi_{2}^{e},\phi_{2}-\varphi\big)_{D^{k}}\\
&+\big(\phi_{1}-\phi_{1}^{e},\phi_{1}-\epsilon\big)_{D^{k}}-\big(\varphi-\varphi^{e},\sigma\big)_{D^{k}} +\big(\vartheta-\vartheta^{e},\tau\big)_{D^{k}}+\lambda_{1}\big(\phi_{1}-\phi_{1}^{e},\pi\big)_{D^{k}}-\lambda_{1}\big(\phi_{2}-\phi_{2}^{e},\sigma\big)_{D^{k}}
\\
&-\big(n.(\tau-\tau^{e})^{*},\pi\big)_{\partial D^{k}}+\big(n.(\sigma-\sigma^{e})^{*},\vartheta\big)_{\partial D^{k}}
+\big(n.(\vartheta-\vartheta^{e})^{*},\sigma\big)_{\partial D^{k}}-\big(n.(\pi-\pi^{e})^{*},\tau\big)_{\partial D^{k}}=0.\\
\end{split}
\end{equation}
Summing over $k$, simplify  by integration by parts and \eqref{fl}. This completes
the proof. $\quad\Box$

\begin{thm}\label{fgrh}
Let $u$ be the exact solution of the problem \eqref{1bchj}, and let $u_{h}$ be the numerical solution of the semi-discrete LDG scheme \eqref{1bch2}.  Then for small enough $h$, we have the following error estimates:
\begin{equation}\label{tt7hb}
\begin{split}
&\|u(.,t)-u_{h}(.,t)\|_{L^{2}(\Omega_{h})}\leq Ch^{N+1},\\
\end{split}
\end{equation}
\end{thm}
where the constant $C$ is dependent upon $T$ and some norms of the solutions. \\
\textbf{Proof}.
Integrating both sides of the above  identity Lemma \ref{lm1h} with respect to $t$ over $(0,T)$, we get
 \begin{equation}\label{tt6h}
\begin{split}
&\frac{1}{2}\|\pi(.,T)\|^{2}_{L^{2}(\Omega_{h})}+\frac{1}{2}\|\sigma(.,T)\|^{2}_{L^{2}(\Omega_{h})}
+ \int_{0}^{T}\big(\big(\Delta_{(\alpha-2)/2}\epsilon,\epsilon\big)_{\Omega_{h}}
+\big(\Delta_{(\alpha-2)/2}\varphi,\varphi\big)_{\Omega_{h}}+\big(\phi_{1},\phi_{1}\big)_{\Omega_{h}}
+\big(\phi_{2},\phi_{2}\big)_{\Omega_{h}}\big)dt\\
&=
\frac{1}{2}\|\pi(.,0)\|^{2}_{L^{2}(\Omega_{h})}
+\frac{1}{2}\|\sigma(.,0)\|^{2}_{L^{2}(\Omega_{h})}+\sum_{j=1}^{4}\int_{0}^{T}Q_{j}dt. \\
\end{split}
\end{equation}
 Next we estimate the term $\int_{0}^{T}Q_{i}dt$, $\,\, i=1,...,4$. So we employ Young's inequality \eqref{91chjm} and the approximation results \eqref{sth}, we obtain
 \begin{equation}\label{tt12h}
\begin{split}
\int_{0}^{T}Q_{1}dt\leq &\int_{0}^{T}(c_{5}\|\epsilon\|^{2}_{L^{2}(\Omega_{h})}
+c_{6} \|\varphi\|^{2}_{L^{2}(\Omega_{h})}+c_{1}\|\pi\|^{2}_{L^{2}(\Omega_{h})}+
c_{2}\|\sigma\|^{2}_{L^{2}(\Omega_{h})}+c_{3}\|\phi_{1}\|^{2}_{L^{2}(\Omega_{h})}
+c_{4}\|\phi_{2}\|^{2}_{L^{2}(\Omega_{h})})dt.\\
\end{split}
\end{equation}
Using the definition of the numerical traces, \eqref{fl}, and the definitions of the projections $\mathcal{P}^{+},\mathcal{P}^{-}$ \eqref{prh}, we get
\begin{equation}\label{91ch}
\begin{split}
Q_{2}=Q_{4}=0.
\end{split}
\end{equation}
So
\begin{equation}\label{tt5h}
\begin{split}
\int_{0}^{T}(Q_{2}+Q_{4})dt=0.
\end{split}
\end{equation}
From the approximation results \eqref{sth} and Young's inequality, we obtain
\begin{equation}\label{tt12h}
\begin{split}
\int_{0}^{T}Q_{3}dt\leq &\int_{0}^{T}(c_{5}\|\epsilon\|^{2}_{L^{2}(\Omega_{h})}
+c_{6} \|\varphi\|^{2}_{L^{2}(\Omega_{h})}+c_{1}\|\pi\|^{2}_{L^{2}(\Omega_{h})}+
c_{2}\|\sigma\|^{2}_{L^{2}(\Omega_{h})})dt\\
&+c_{3}\|\phi_{1}\|^{2}_{L^{2}(\Omega_{h})}
+c_{4}\|\phi_{2}\|^{2}_{L^{2}(\Omega_{h})}+Ch^{2N+2}.\\
\end{split}
\end{equation}

Combining \eqref{tt12h}, \eqref{tt5h} and \eqref{tt6h}, we obtain
\begin{equation}\label{91ch}
\begin{split}
&\frac{1}{2}\|\pi(.,T)\|^{2}_{L^{2}(\Omega_{h})}+\frac{1}{2}\|\sigma(.,T)\|^{2}_{L^{2}(\Omega_{h})}
+ \int_{0}^{T}\big(\big(\Delta_{(\alpha-2)/2}\epsilon,\epsilon\big)_{\Omega_{h}}
+\big(\Delta_{(\alpha-2)/2}\varphi,\varphi\big)_{\Omega_{h}}+\big(\phi_{1},\phi_{1}\big)_{\Omega_{h}}
+\big(\phi_{2},\phi_{2}\big)_{\Omega_{h}}\big)dt\\
&\leq \frac{1}{2}\|\pi(.,0)\|^{2}_{L^{2}(\Omega_{h})}
+\frac{1}{2}\|\sigma(.,0)\|^{2}_{L^{2}(\Omega_{h})}+\int_{0}^{T}(c_{1}\|\pi\|^{2}_{L^{2}(\Omega_{h})}
+c_{2}\|\sigma\|^{2}_{L^{2}(\Omega_{h})})dt+\int_{0}^{T}(c_{5}\|\epsilon\|^{2}_{L^{2}(\Omega_{h})}\\
&
+c_{6} \|\varphi\|^{2}_{L^{2}(\Omega_{h})}+c_{3}\|\phi_{1}\|^{2}_{L^{2}(\Omega_{h})}
+c_{4}\|\phi_{2}\|^{2}_{L^{2}(\Omega_{h})})dt
+Ch^{2N+2}. \\
\end{split}
\end{equation}
Recalling Lemmas \ref{lg}, we obtain
\begin{equation}\label{91ch}
\begin{split}
&\frac{1}{2}\|\pi(.,T)\|^{2}_{L^{2}(\Omega_{h})}+\frac{1}{2}\|\sigma(.,T)\|^{2}_{L^{2}(\Omega_{h})}
+ \int_{0}^{T}\big(\big(\phi_{1},\phi_{1}\big)_{\Omega_{h}}
+\big(\phi_{2},\phi_{2}\big)_{\Omega_{h}}\big)dt\\
&\leq \frac{1}{2}\|\pi(.,0)\|^{2}_{L^{2}(\Omega_{h})}
+\frac{1}{2}\|\sigma(.,0)\|^{2}_{L^{2}(\Omega_{h})}+\int_{0}^{T}(c_{1}\|\pi\|^{2}_{L^{2}(\Omega_{h})}
+c_{2}\|\sigma\|^{2}_{L^{2}(\Omega_{h})})dt\\
&\quad\quad\int_{0}^{T}(c_{3}\|\phi_{1}\|^{2}_{L^{2}(\Omega_{h})}
+c_{4}\|\phi_{2}\|^{2}_{L^{2}(\Omega_{h})})dt
+Ch^{2N+2}, \\
\end{split}
\end{equation}
provided $c_{i},\,\,i=1,2,3,4$ are sufficiently small such that $c_{i}\leq1$, we obtain
\begin{equation}\label{91ch}
\begin{split}
&\frac{1}{2}\|\pi(.,T)\|^{2}_{L^{2}(\Omega_{h})}+\frac{1}{2}\|\sigma(.,T)\|^{2}_{L^{2}(\Omega_{h})}
\\
&\leq \frac{1}{2}\|\pi(.,0)\|^{2}_{L^{2}(\Omega_{h})}
+\frac{1}{2}\|\sigma(.,0)\|^{2}_{L^{2}(\Omega_{h})}+\int_{0}^{T}(\|\pi\|^{2}_{L^{2}(\Omega_{h})}
+\|\sigma\|^{2}_{L^{2}(\Omega_{h})})dt
+Ch^{2N+2}. \\
\end{split}
\end{equation}
 Employing Gronwall's lemma, we can get \eqref{tt7hb}. $\quad\Box$

\section{ LDG method for  strongly nonlinear  coupled
fractional Schr\"{o}dinger equations}\label{sc40}

In this section, we present and analyze the LDG method for the strongly  coupled  nonlinear fractional Schr\"{o}dinger equations

\begin{equation}\label{1bchcc}
\begin{split}
&i\frac{\partial u_{1}}{\partial t}- \lambda_{1}(-\Delta)^{\frac{\alpha}{2}}u_{1}+ \varpi_{1}u_{1}+ \varpi_{2}u_{2}+ \lambda_{2}f(|u_{1}|^{2},|u_{2}|^{2})u_{1}=0,\\
&i\frac{\partial u_{2}}{\partial t}- \lambda_{3}(-\Delta)^{\frac{\alpha}{2}}u_{2}+\varpi_{2} u_{1}+ \varpi_{1}u_{2}+\lambda_{4} g(|u_{1}|^{2},|u_{2}|^{2})u_{2}=0.\\
\end{split}
\end{equation}
To define the local discontinuous Galerkin method, we rewrite  \eqref{1bchcc} as a first-order system:

\begin{equation}\label{1bch}
\begin{split}
&i\frac{\partial u_{1}}{\partial t}+\lambda_{1}e+ \varpi_{1}u_{1}+\varpi_{2} u_{2}+ \lambda_{2}f(|u_{1}|^{2},|u_{2}|^{2})u_{1}=0,\\
&e=\Delta_{(\alpha-2)/2}r,\quad r=\frac{\partial}{\partial x}s,\quad s=\frac{\partial}{\partial x}u_{1},\\
&i\frac{\partial u_{2}}{\partial t}+\lambda_{3}l+ \varpi_{2}u_{1}+\varpi_{1} u_{2}+ \lambda_{4}g(|u_{1}|^{2},|u_{2}|^{2})u_{2}=0,\\
&l=\Delta_{(\alpha-2)/2}w,\quad w=\frac{\partial}{\partial x}z,\quad z=\frac{\partial}{\partial x}u_{2}.\\
\end{split}
\end{equation}
We decompose the complex functions $u(x, t)$ and $v(x, t)$ into their real and imaginary parts. Setting
$u_{1}(x, t)=p(x, t) + iq(x, t)$ and $u_{2}(x, t) = \upsilon(x, t) + i\theta(x, t)$ in system \eqref{1bchcc}, we can obtain the following coupled system
\begin{equation}\label{1bx}
\begin{split}
&\frac{\partial p}{\partial t}+\lambda_{1}e_{1}+ \varpi_{1}q+ \varpi_{2}\theta+ \lambda_{2}f(|u_{1}|^{2},|u_{2}|^{2})q=0,\\
&e_{1}=\Delta_{(\alpha-2)/2}r,\quad r=\frac{\partial}{\partial x}s,\quad s=\frac{\partial}{\partial x}q,\\
&\frac{\partial q}{\partial t}-\lambda_{1}l_{1}- \varpi_{1}p- \varpi_{2}\upsilon- \lambda_{2}f(|u_{1}|^{2},|u_{2}|^{2})p=0,\\
&l_{1}=\Delta_{(\alpha-2)/2}w,\quad w=\frac{\partial}{\partial x}z,\quad z=\frac{\partial}{\partial x}p,\\
&\frac{\partial \upsilon}{\partial t}+\lambda_{3}e_{2}+\varpi_{3}q+ \varpi_{4}\theta+ \lambda_{4}g(|u_{1}|^{2},|u_{2}|^{2})\theta=0,\\
&e_{2}=\Delta_{(\alpha-2)/2}\rho,\quad\rho=\frac{\partial}{\partial x}\varpi,\quad \varpi=\frac{\partial}{\partial x}\theta,\\
&\frac{\partial \theta}{\partial t}-\lambda_{3}l_{2}-\varpi_{2} p- \varpi_{1}\upsilon- \lambda_{4}g(|u_{1}|^{2},|u_{2}|^{2})\upsilon=0,\\
&l_{2}=\Delta_{(\alpha-2)/2}\xi,\quad\xi=\frac{\partial}{\partial x}\varrho,\quad \varrho=\frac{\partial}{\partial x}\upsilon.\\
\end{split}
\end{equation}

We define  local discontinuous Galerkin scheme as follows: find $p_{h}, q_{h}, e_{1},r_{h},s_{h},l_{1},w_{h},z_{h}$,\\
$\upsilon_{h},\theta_{h},e_{2},\rho_{h},\varpi_{h},l_{2},\xi_{h}$,$\varrho_{h} \in V_{k}^{N}$, such that for all test functions $\vartheta_{1},\beta_{1},\phi,\varphi,\chi,\beta_{2},\psi$,
$\zeta,\gamma,\beta_{3},\delta,\varsigma,o,\beta_{4},\omega,\kappa\in V_{k}^{N}$,
\begin{equation*}\label{1bxhh}
\begin{split}
&\big(\frac{\partial p_{h}}{\partial t},\vartheta_{1}\big)_{D^{k}}+\lambda_{1}\big(T_{h},\vartheta_{1}\big)_{D^{k}}+ \varpi_{1}\big(q_{h},\vartheta_{1}\big)_{D^{k}}+ \varpi_{2}\big(\theta_{h},\vartheta_{1}\big)_{D^{k}}+ \lambda_{2}\big(f(|u_{1}|^{2},|u_{2}|^{2})q_{h},\vartheta_{1}\big)_{D^{k}}=0,\\
&\big(T_{h},\beta_{1}\big)_{D^{k}}=\big(\Delta_{(\alpha-2)/2}r_{h},\beta_{1}\big)_{D^{k}},\\
&\big(r_{h},\phi\big)_{D^{k}}=\big(\frac{\partial}{\partial x}s_{h},\phi\big)_{D^{k}},\\
&\big(s_{h},\varphi\big)_{D^{k}}=\big(\frac{\partial}{\partial x}q_{h},\varphi\big)_{D^{k}},\\
&\big(\frac{\partial q_{h}}{\partial t},\chi\big)_{D^{k}}-\lambda_{1}\big(H_{h},\chi\big)_{D^{k}}
- \varpi_{1}\big(p_{h},\chi\big)_{D^{k}}-\varpi_{2} \big(\upsilon_{h},\chi\big)_{D^{k}}- \lambda_{2}\big(f(|u_{1}|^{2},|u_{2}|^{2})p_{h},v\big)_{D^{k}}=0,\\
&\big(H_{h},\beta_{2}\big)_{D^{k}}=\big(\Delta_{(\alpha-2)/2}w_{h},\beta_{2}\big)_{D^{k}},\\
&\big(w_{h},\psi\big)_{D^{k}}=\big(\frac{\partial}{\partial x}z_{h},\psi\big)_{D^{k}},\\
&\big(z_{h},\zeta\big)_{D^{k}}=\big(\frac{\partial}{\partial x}p_{h},\zeta\big)_{D^{k}},\\
\end{split}
\end{equation*}
\begin{equation}\label{1bxhh23}
\begin{split}
&\big(\frac{\partial \upsilon_{h}}{\partial t},\gamma\big)_{D^{k}}+\lambda_{3}\big(L_{h},\gamma\big)_{D^{k}}
+\varpi_{2}\big(q_{h},\gamma\big)_{D^{k}}+ \varpi_{1}\big(\theta_{h},\gamma\big)_{D^{k}}+ \lambda_{4}\big(g(|u_{1}|^{2},|u_{2}|^{2})\theta_{h},\gamma\big)_{D^{k}}=0,\\
&\big(L_{h},\beta_{3}\big)_{D^{k}}=\big(\Delta_{(\alpha-2)/2}\rho_{h},\beta_{3}\big)_{D^{k}},\\
&\big(\rho_{h},\delta\big)_{D^{k}}=\big(\frac{\partial}{\partial x}\varpi_{h},\delta\big)_{D^{k}},\\
&\big(\varpi_{h},\varsigma\big)_{D^{k}}=\big(\frac{\partial}{\partial x}\theta_{h},\varsigma\big)_{D^{k}},\\
&\big(\frac{\partial \theta_{h}}{\partial t},o\big)_{D^{k}}-\lambda_{3}\big(E_{h},o\big)_{D^{k}}-\varpi_{2}\big(p_{h},o\big)_{D^{k}}- \varpi_{1}\big(\upsilon_{h},o\big)_{D^{k}}- \lambda_{4}\big(g(|u_{1}|^{2},|u_{2}|^{2})\upsilon_{h},o\big)_{D^{k}}=0,\\
&\big(E_{h},\beta_{4}\big)_{D^{k}}=\big(\Delta_{(\alpha-2)/2}\xi_{h},\beta_{4}\big)_{D^{k}},\\
&\big(\xi_{h},\omega\big)_{D^{k}}=\big(\frac{\partial}{\partial x}\varrho_{h},\omega\big)_{D^{k}},\\
&\big(\varrho_{h},\kappa\big)_{D^{k}}=\big(\frac{\partial}{\partial x}\upsilon_{h},\kappa\big)_{D^{k}}.\\
\end{split}
\end{equation}
Applying integration by parts to \eqref{1bxhh23}, and replacing the fluxes at the interfaces
by the corresponding numerical fluxes, we obtain
\begin{equation}\label{1bch2cc}
\begin{split}
&\big(\frac{\partial p_{h}}{\partial t},\vartheta_{1}\big)_{D^{k}}+\lambda_{1}\big(T_{h},\vartheta_{1}\big)_{D^{k}}+ \varpi_{1}\big(q_{h},\vartheta_{1}\big)_{D^{k}}+\varpi_{2}\big(\theta_{h},\vartheta_{1}\big)_{D^{k}}+ \lambda_{2}\big(f(|u_{1}|^{2},|u_{2}|^{2})q_{h},\vartheta_{1}\big)_{D^{k}}=0,\\
&\big(T_{h},\beta_{1}\big)_{D^{k}}=\big(\Delta_{(\alpha-2)/2}r_{h},\beta_{1}\big)_{D^{k}},\\
&\big(r_{h},\phi\big)_{D^{k}}=-\big(s_{h},\phi_{x}\big)_{D^{k}}+\big(n.s^{*}_{h},\phi\big)_{\partial D^{k}},\\
&\big(s_{h},\varphi\big)_{D^{k}}=-\big(q_{h},\varphi_{x}\big)_{D^{k}}+\big(n.q^{*}_{h},\varphi\big)_{\partial D^{k}},\\
&\big(\frac{\partial q_{h}}{\partial t},\chi\big)_{D^{k}}-\lambda_{1}\big(H_{h},\chi\big)_{D^{k}}
-\varpi_{1} \big(p_{h},\chi\big)_{D^{k}}- \varpi_{2}\big(\upsilon_{h},\chi\big)_{D^{k}}- \lambda_{2}\big(f(|u_{1}|^{2},|u_{2}|^{2})p_{h},\chi\big)_{D^{k}}=0,\\
&\big(H_{h},\beta_{2}\big)_{D^{k}}=\big(\Delta_{(\alpha-2)/2}w_{h},\beta_{2}\big)_{D^{k}},\\
&\big(w_{h},\psi\big)_{D^{k}}=-\big(z_{h},\psi_{x}\big)_{D^{k}}+\big(n.z^{*}_{h},\psi\big)_{\partial D^{k}},\\
&\big(z_{h},\zeta\big)_{D^{k}}=-\big(p_{h},\zeta_{x}\big)_{D^{k}}+\big(n.p^{*}_{h},\zeta\big)_{\partial D^{k}},\\
&\big(\frac{\partial \upsilon_{h}}{\partial t},\gamma\big)_{D^{k}}+\lambda_{3}\big(L_{h},\gamma\big)_{D^{k}}
+\varpi_{2}\big(q_{h},\gamma\big)_{D^{k}}+\varpi_{1}\big(\theta_{h},\gamma\big)_{D^{k}}+ \lambda_{4}\big(g(|u_{1}|^{2},|u_{2}|^{2})\theta_{h},\gamma\big)_{D^{k}}=0,\\
&\big(L_{h},\beta_{3}\big)_{D^{k}}=\big(\Delta_{(\alpha-2)/2}\rho_{h},\beta_{3}\big)_{D^{k}},\\
&\big(\rho_{h},\delta\big)_{D^{k}}=-\big(\varpi_{h},\delta_{x}\big)_{D^{k}}+\big(n.\varpi_{h}^{*},\delta\big)_{\partial D^{k}},\\
&\big(\varpi_{h},\varsigma\big)_{D^{k}}=-\big(\theta_{h},\varsigma_{x}\big)_{D^{k}}+\big(n.\theta_{h}^{*},\varsigma\big)_{\partial D^{k}},\\
&\big(\frac{\partial \theta_{h}}{\partial t},o\big)_{D^{k}}-\lambda_{3}\big(E_{h},o\big)_{D^{k}}-\varpi_{2}\big(p_{h},o\big)_{D^{k}}- \varpi_{1} \big(\upsilon_{h},o\big)_{D^{k}}- \lambda_{4}\big(g(|u_{1}|^{2},|u_{2}|^{2})\upsilon_{h},o\big)_{D^{k}}=0,\\
&\big(E_{h},\beta_{4}\big)_{D^{k}}=\big(\Delta_{(\alpha-2)/2}\xi_{h},\beta_{4}\big)_{D^{k}},\\
&\big(\xi_{h},\omega\big)_{D^{k}}=-\big(\varrho_{h},\omega_{x}\big)_{D^{k}}
+\big(n.\varrho_{h}^{*},\omega\big)_{\partial D^{k}},\\
&\big(\varrho_{h},\kappa\big)_{D^{k}}=-\big(\upsilon_{h},\kappa_{x}\big)_{D^{k}}
+\big(n.\upsilon_{h}^{*},\kappa\big)_{\partial D^{k}},\\
\end{split}
\end{equation}
The numerical traces $( p, q, s,z,\upsilon,\theta,\varpi,\varrho)$ are defined on interelement faces as the
alternating fluxes
\begin{equation}\label{flc}
\begin{split}
&p^{*}_{k+\frac{1}{2}}=p^{-}_{k+\frac{1}{2}},\,\,s^{*}_{k+\frac{1}{2}}
=s^{+}_{k+\frac{1}{2}},\,\,\,\,
q^{*}_{k+\frac{1}{2}}=q^{-}_{k+\frac{1}{2}},
z^{*}_{k+\frac{1}{2}}=z^{+}_{k+\frac{1}{2}},\\
& \upsilon^{*}_{k+\frac{1}{2}}=\upsilon^{-}_{k+\frac{1}{2}},\,\,\varpi^{*}_{k+\frac{1}{2}}
=\varpi^{+}_{k+\frac{1}{2}},\,\,\,\,\varrho^{*}_{k+\frac{1}{2}}=\varrho^{+}_{k+\frac{1}{2}},\,\,
\theta^{*}_{k+\frac{1}{2}}=\theta^{-}_{k+\frac{1}{2}}.
\end{split}
\end{equation}

\section{ Stability and error estimates}\label{sc4}
 In the following we discuss stability and accuracy of the proposed scheme, for the nonlinear fractional coupled
Schr\"{o}dinger  problem.
\subsection{ Stability analysis }  In order to carry out the analysis of the LDG scheme,
\begin{thm}\label{tt4h} ($L^{2}$ stability).
The semidiscrete scheme \eqref{1bch2cc} is   stable, and \\
$
\|u_{h}(x,T)\|_{\Omega_{h}}+\|v_{h}(x,T)\|_{\Omega_{h}}\leq c(\|u_{0}(x)\|_{\Omega_{h}}+\|v_{0}(x)\|_{\Omega_{h}})$ for any $T>0$.
\end{thm}
\textbf{Proof.}
Set $(\vartheta_{1},\beta_{1},\phi,\varphi,\chi,\beta_{2},\psi,\zeta,\gamma,\beta_{3},\delta,\varsigma,\beta_{4},o,\omega,\kappa)
=(p_{h},T_{h}-r_{h},p_{h},-z_{h},q_{h},H_{h}-w_{h},-q_{h},s_{h},\upsilon_{h},L_{h}-\rho_{h},
\upsilon_{h},-\varrho_{h},\theta_{h},E_{h}-\xi_{h},-\theta_{h},\varpi_{h})$ in \eqref{1bch2},
 and consider the integration by parts formula $\big(u,\frac{\partial r}{\partial x}\big)_{D^{k}}+\big(r,\frac{\partial u}{\partial x}\big)_{D^{k}}=[ur]_{_{x_{k-\frac{1}{2}}}}^{x_{k+\frac{1}{2}}}$, we get
\begin{equation}\label{1bch1bn}
\begin{split}
&\big((p_{h})_{t},p_{h}\big)_{D^{k}}+\big((q_{h})_{t},q_{h}\big)_{D^{k}}
+\big((\upsilon_{h})_{t},\upsilon_{h}\big)_{D^{k}}+\big((\theta_{h})_{t},\theta_{h}\big)_{D^{k}}
+\big(\Delta_{(\alpha-2)/2}w_{h},w_{h}\big)_{D^{k}}+\big(\Delta_{(\alpha-2)/2}\xi_{h},\xi_{h}\big)_{D^{k}}\\
&\quad+\big(\Delta_{(\alpha-2)/2}r_{h},r_{h}\big)_{D^{k}}
+\big(\Delta_{(\alpha-2)/2}\rho_{h},\rho_{h}\big)_{D^{k}}+\big(T_{h},T_{h}\big)_{D^{k}}+\big(H_{h},H_{h}\big)_{D^{k}}
+\big(L_{h},L_{h}\big)_{D^{k}}+\big(E_{h},E_{h}\big)_{D^{k}}\\
&=
\big(\Delta_{(\alpha-2)/2}w_{h},H_{h}\big)_{D^{k}}+\big(\Delta_{(\alpha-2)/2}\xi_{h},E_{h}\big)_{D^{k}}
+\big(\Delta_{(\alpha-2)/2}r_{h},T_{h}\big)_{D^{k}}
+\big(\Delta_{(\alpha-2)/2}\rho_{h},L_{h}\big)_{D^{k}}\\
&\quad-\big( T_{h},-r_{h}+\lambda_{1}p_{h}\big)_{D^{k}}+\big( H_{h},w_{h}+\lambda_{1}q_{h}\big)_{D^{k}}-\big( L_{h},\lambda_{3}\upsilon_{h}-\rho_{h}\big)_{D^{k}}
+\big( E_{h},\xi_{h}+\lambda_{3}\theta_{h}\big)_{D^{k}}\\
&\quad-\big( r_{h},p_{h}\big)_{D^{k}}-\big( \rho_{h},\upsilon_{h}\big)_{D^{k}}+\big( w_{h},q_{h}\big)_{D^{k}}+\big( \xi_{h},\theta_{h}\big)_{D^{k}}+\theta(s_{h},p_{h})+\theta(\varpi_{h},\upsilon_{h})-\theta(q_{h},z_{h})
-\theta(\theta_{h},\varrho_{h}).\\
\end{split}
\end{equation}

Summing over all elements \eqref{1bch1bn}, employing Young's inequality and using the definition of the numerical traces, \eqref{flc},  we obtain
\begin{equation}\label{1bch1}
\begin{split}
&\big((p_{h})_{t},p_{h}\big)_{\Omega_{h}}+\big((q_{h})_{t},q_{h}\big)_{\Omega_{h}}
+\big((\upsilon_{h})_{t},\upsilon_{h}\big)_{\Omega_{h}}+\big((\theta_{h})_{t},\theta_{h}\big)_{\Omega_{h}}
+\big(\Delta_{(\alpha-2)/2}w_{h},w_{h}\big)_{\Omega_{h}}+\big(\Delta_{(\alpha-2)/2}\xi_{h},\xi_{h}\big)_{\Omega_{h}}\\
&\quad+\big(\Delta_{(\alpha-2)/2}r_{h},r_{h}\big)_{\Omega_{h}}+\big(\Delta_{(\alpha-2)/2}\rho_{h},\rho_{h}\big)_{D^{k}}+\big(H_{h},H_{h}\big)_{\Omega_{h}}
+\big(L_{h},L_{h}\big)_{\Omega_{h}}+\big(E_{h},E_{h}\big)_{\Omega_{h}}+\big(T_{h},T_{h}\big)_{\Omega_{h}}\\
&\leq
c_{12}\|w_{h}\|^{2}_{L^{2}(\Omega_{h})}+
c_{11}\|r_{h}\|^{2}_{L^{2}(\Omega_{h})}+c_{10}\|\xi_{h}\|^{2}_{L^{2}(\Omega_{h})}+
c_{9}\|\rho_{h}\|^{2}_{L^{2}(\Omega_{h})}+c_{5}\|p_{h}\|^{2}_{L^{2}(\Omega_{h})}
+c_{6}\|q_{h}\|^{2}_{L^{2}(\Omega_{h})}
+c_{7}\|\upsilon_{h}\|^{2}_{L^{2}(\Omega_{h})}\\
&\quad+c_{8}\|\theta_{h}\|^{2}_{L^{2}(\Omega_{h})}+c_{1}\|T_{h}\|^{2}_{L^{2}(\Omega_{h})}+c_{2}\|H_{h}\|^{2}_{L^{2}(\Omega_{h})}
+c_{3}\|E_{h}\|^{2}_{L^{2}(\Omega_{h})}
+c_{4}\|L_{h}\|^{2}_{L^{2}(\Omega_{h})}.
\end{split}
\end{equation}
Recalling Lemma \ref{lg} and provided $c_{i},\,\,i=1,2,...,8$ are sufficiently small such that $c_{i}\leq1$, we obtain that

\begin{equation}\label{1bch3}
\begin{split}
\big((p_{h})_{t},p_{h}\big)_{\Omega_{h}}+\big((q_{h})_{t},q_{h}\big)_{\Omega_{h}}
+\big((\upsilon_{h})_{t},\upsilon_{h}\big)_{\Omega_{h}}+\big((\theta_{h})_{t},\theta_{h}\big)_{\Omega_{h}}\leq
&\|p_{h}\|^{2}_{L^{2}(\Omega_{h})}
+\|q_{h}\|^{2}_{L^{2}(\Omega_{h})}\\
&+\|\upsilon_{h}\|^{2}_{L^{2}(\Omega_{h})}+\|\theta_{h}\|^{2}_{L^{2}(\Omega_{h})}.
\end{split}
\end{equation}
Hence
\begin{equation}\label{1bch3}
\begin{split}
\frac{1}{2}\frac{d}{dt}\|u_{h}\|^{2}_{\Omega_{h}}+\frac{1}{2}\frac{d}{dt}\|v_{h}\|^{2}_{\Omega_{h}}\leq
&\|u\|^{2}_{\Omega_{h}}+\|v\|^{2}_{\Omega_{h}}.\\
\end{split}
\end{equation}
Employing Gronwall's inequality, we obtain
\begin{equation}\label{cd}
\begin{split}
\|u_{h}(x,T)\|_{\Omega_{h}}^{2}+\|v_{h}(x,T)\|_{\Omega_{h}}^{2}\leq C((\|u_{0}(x)\|_{\Omega_{h}}^{2}+\|v_{0}(x)\|_{\Omega_{h}}^{2}). \\
\end{split}
\end{equation}

\subsection{Error estimates}
We consider the linear  fractional coupled Schr\"{o}dinger system
\begin{equation}\label{1bchjcn}
\begin{split}
&i\frac{\partial u_{1}}{\partial t}-\lambda_{1}(-\Delta)^{\frac{\alpha}{2}}u_{1}+\omega_{1}u_{1}+\omega_{2}u_{2}+\lambda_{2}u_{1}=0,\\
&i\frac{\partial u_{2}}{\partial t}-\lambda_{3}(-\Delta)^{\frac{\alpha}{2}}u_{2}+\omega_{2}u_{1}+\omega_{1}u_{2}+\lambda_{4}u_{2}=0.\\
\end{split}
\end{equation}
It is easy to verify that the error equations of the above  \eqref{1bchjcn} satisfies

\begin{equation}\label{1bxhhnm}
\begin{split}
\big(&\frac{\partial (p-p_{h})}{\partial t},\vartheta_{1}\big)_{D^{k}}+\big(\frac{\partial (q-q_{h})}{\partial t},\chi\big)_{D^{k}}+\big(\frac{\partial (\upsilon-\upsilon_{h})}{\partial t},\gamma\big)_{D^{k}}+\big(\frac{\partial (\theta-\theta_{h})}{\partial t},o\big)_{D^{k}}-\big(\Delta_{(\alpha-2)/2}(r-r_{h}),\beta_{1}\big)_{D^{k}}\\
&-\big(\Delta_{(\alpha-2)/2}(w-w_{h}),\beta_{2}\big)_{D^{k}}-
\big(\Delta_{(\alpha-2)/2}(\rho-\rho_{h}),\beta_{3}\big)_{D^{k}}
-\big(\Delta_{(\alpha-2)/2}(\xi-\xi_{h}),\beta_{4}\big)_{D^{k}}
\\
&+\lambda_{1}\big(T-T_{h},\vartheta_{1}\big)_{D^{k}}-\lambda_{1}\big(H-H_{h},\chi\big)_{D^{k}}+
\lambda_{3}\big(L-L_{h},\gamma\big)_{D^{k}}-\lambda_{3}\big(E-E_{h},o\big)_{D^{k}}
+\big(T-T_{h},\beta_{1}\big)_{D^{k}}\\
&+\big(H-H_{h},\beta_{2}\big)_{D^{k}}+
\big(L-L_{h},\beta_{3}\big)_{D^{k}}+\big(E-E_{h},\beta_{4}\big)_{D^{k}}
+\big(q-q_{h},\varphi_{x}\big)_{D^{k}}+\big(s-s_{h},\phi_{x}\big)_{D^{k}}\\
&+\big(z-z_{h},\psi_{x}\big)_{D^{k}}
+\big(p-p_{h},\zeta_{x}\big)_{D^{k}}+\big(\varpi-\varpi_{h},\delta_{x}\big)_{D^{k}}
+\big(\theta-\theta_{h},\varsigma_{x}\big)_{D^{k}}
+\big(\upsilon-\upsilon_{h},\kappa_{x}\big)_{D^{k}}\\
&+\big(\varrho-\varrho_{h},\omega_{x}\big)_{D^{k}}+ \omega_{1}\big(q-q_{h},\vartheta_{1}\big)_{D^{k}}+ \omega_{2}\big(\theta-\theta_{h},\vartheta_{1}\big)_{D^{k}}+ \lambda_{2}\big(q-q_{h},\vartheta_{1}\big)_{D^{k}}+\big(r-r_{h},\phi\big)_{D^{k}}+\big(s-s_{h},\varphi\big)_{D^{k}}
\\
&- \omega_{1}\big(p-p_{h},\chi\big)_{D^{k}}- \omega_{2}\big(\upsilon-\upsilon_{h},\chi\big)_{D^{k}}- \lambda_{2}\big(p-p_{h},\chi\big)_{D^{k}}+\big(w-w_{h},\psi\big)_{D^{k}}+\big(z-z_{h},\zeta\big)_{D^{k}}
+\omega_{2}\big(q-q_{h},\gamma\big)_{D^{k}}\\
&+ \omega_{1}\big(\theta-\theta_{h},\gamma\big)_{D^{k}}+ \lambda_{4} \big(\theta-\theta_{h},\gamma\big)_{D^{k}}+\big(\rho-\rho_{h},\delta\big)_{D^{k}}
+\big(\varpi-\varpi_{h},\varsigma\big)_{D^{k}}-\omega_{2}\big(p-p_{h},o\big)_{D^{k}}- \omega_{1}\big(\upsilon-\upsilon_{h},o\big)_{D^{k}}\\
&- \lambda_{4}\big(\upsilon-\upsilon_{h},o\big)_{D^{k}}+\big(\xi-\xi_{h},\omega\big)_{D^{k}}
+\big(\varrho-\varrho_{h},\kappa\big)_{D^{k}}
-\big(n.(s-s_{h})^{*},\phi\big)_{\partial D^{k}}-\big(n.(q-q_{h})^{*},\varphi\big)_{\partial D^{k}}\\
&-\big(n.(z-z_{h})^{*},\psi\big)_{\partial D^{k}}-\big(n.(p-p_{h})^{*},\zeta\big)_{\partial D^{k}}-\big(n.(\varpi-\varpi_{h})^{*},\delta\big)_{\partial D^{k}}-\big(n.(\theta-\theta_{h})^{*},\varsigma\big)_{\partial D^{k}}\\
&+\big(n.(\varrho-\varrho_{h})^{*},\omega\big)_{\partial D^{k}}
-\big(n.(\upsilon-\upsilon_{h})^{*},\kappa\big)_{\partial D^{k}}=0.\\
\end{split}
\end{equation}

\begin{thm}\label{fgrh}
Let $u$ and $v$  be the exact solutions of the linear coupled  fractional
Schr\"{o}dinger equations \eqref{1bchjcn}, and let $u_{h}$ and $v_{h}$ be the numerical solutions of the semi-discrete LDG scheme \eqref{1bch2cc}.  Then for small enough $h$, we have the following error estimates:
\begin{equation}\label{tt7h}
\begin{split}
&\|u(.,T)-u_{h}(.,T)\|_{L^{2}(\Omega_{h})}+\|v(.,T)-v_{h}(.,T)\|_{L^{2}(\Omega_{h})}\leq Ch^{N+1},\\
\end{split}
\end{equation}
\end{thm}
where the constant $C$ is dependent upon $T$ and some norms of the solutions. \\
\textbf{Proof.}
We donate

\begin{equation}\label{91h}
\begin{split}
&\pi_{1}=\mathcal{P}^{-}\upsilon-\upsilon_{h},\quad \pi^{e}_{1}=\mathcal{P}^{-}\upsilon-\upsilon,\quad \pi_{2}=\mathcal{P}^{-}\theta-\theta_{h},\quad \pi_{2}^{e}=\mathcal{P}^{-}\theta-\theta ,\\
& \pi_{3}=\mathcal{P}^{+}\rho-\rho_{h},\quad\pi_{3}^{e}=\mathcal{P}^{+}\rho-\rho,\quad \pi_{4}=\mathcal{P}^{+}\varpi-\varpi_{h}, \quad \pi_{4}^{e}=\mathcal{P}^{+}\varpi-\varpi,  \\
&\pi_{5}=\mathcal{P}^{+}\xi-\xi_{h},\quad \pi_{5}^{e}=\mathcal{P}^{+}\xi-\xi,\quad \pi_{6}=\mathcal{P}^{+}\varrho-\varrho_{h},\quad \pi_{6}^{e}=\mathcal{P}^{+}\varrho-\varrho,\\
&\epsilon_{1}=\mathcal{P}^{+}T-T_{h},\quad \epsilon_{1}^{e}=\mathcal{P}^{+}T-T,\quad\epsilon_{3}=\mathcal{P}^{+}H-H_{h},
\quad\epsilon_{3}^{e}=\mathcal{P}^{+}H-H,\\
& \epsilon_{4}=\mathcal{P}^{+}L-L_{h}, \quad \epsilon_{4}^{e}=\mathcal{P}^{+}L-L,  \quad \epsilon_{5}=\mathcal{P}^{+}E-E_{h},\quad \epsilon_{5}^{e}=\mathcal{P}^{+}E-E. 
\end{split}
\end{equation}
From the Galerkin orthogonality \eqref{1bxhhnm}, we get
\begin{equation}\label{1bxhhn}
\begin{split}
\big(&\frac{\partial (\pi-\pi^{e})}{\partial t},\vartheta_{1}\big)_{D^{k}}+\big(\frac{\partial (\sigma-\sigma^{e})}{\partial t},\chi\big)_{D^{k}}+\big(\frac{\partial (\pi_{1}-\pi_{1}^{e})}{\partial t},\gamma\big)_{D^{k}}+\big(\frac{\partial (\pi_{2}-\pi_{2}^{e})}{\partial t},o\big)_{D^{k}}
-\big(\Delta_{(\alpha-2)/2}(\epsilon-\epsilon^{e}),\beta_{1}\big)_{D^{k}}\\
&-\big(\Delta_{(\alpha-2)/2}(\varphi-\varphi^{e}),\beta_{2}\big)_{D^{k}}-
\big(\Delta_{(\alpha-2)/2}(\pi_{3}-\pi_{3}^{e}),\beta_{3}\big)_{D^{k}}
-\big(\Delta_{(\alpha-2)/2}(\pi_{5}-\pi_{5}^{e}),\beta_{4}\big)_{D^{k}}
\\
&+\lambda_{1}\big(\epsilon_{1}-\epsilon_{1}^{e},\vartheta\big)_{D^{k}}-
\lambda_{1}\big(\epsilon_{2}-\epsilon_{2}^{e},\chi\big)_{D^{k}}+
\lambda_{3}\big(\epsilon_{3}-\epsilon_{3}^{e},\gamma\big)_{D^{k}}-\lambda_{3}\big(\epsilon_{4}-\epsilon_{4}^{e},o\big)_{D^{k}}
+\big(\epsilon_{1}-\epsilon_{1}^{e},\beta_{1}\big)_{D^{k}}\\
&+\big(\epsilon_{2}-\epsilon_{2}^{e},\beta_{2}\big)_{D^{k}}
+\big(\epsilon_{3}-\epsilon_{3}^{e},\beta_{3}\big)_{D^{k}}+\big(\epsilon_{4}-\epsilon_{4}^{e},\beta_{4}\big)_{D^{k}}
+\big(\tau-\tau^{h},\phi_{x}\big)_{D^{k}}+\big(\sigma-\sigma^{h},\varphi_{x}\big)_{D^{k}}\\
&+\big(\vartheta-\vartheta^{e},\psi_{x}\big)_{D^{k}}
+\big(\pi-\pi^{e},\zeta_{x}\big)_{D^{k}}+\big(\pi_{4}-\pi_{4}^{e},\delta_{x}\big)_{D^{k}}
+\big(\pi_{2}-\pi_{2}^{e},\varsigma_{x}\big)_{D^{k}}
+\big(\pi_{1}-\pi_{1}^{e},\kappa_{x}\big)_{D^{k}}\\
&+\big(\pi_{6}-\pi_{6}^{e},\omega_{x}\big)_{D^{k}}+ \omega_{1}\big(\sigma-\sigma^{e},\vartheta_{1}\big)_{D^{k}}+ \omega_{2}\big(\pi_{2}-\pi_{2}^{e},\vartheta_{1}\big)_{D^{k}}+ \lambda_{2}\big(\sigma-\sigma^{e},\vartheta_{1}\big)_{D^{k}}+\big(\epsilon-\epsilon^{e},\phi\big)_{D^{k}}
+\big(\tau-\tau^{e},\varphi\big)_{D^{k}}
\\
&- \lambda_{2}\big(\pi-\pi^{e},\chi\big)_{D^{k}}-\omega_{1}\big(\pi_{1}-\pi_{1}^{e},\chi\big)_{D^{k}}- \omega_{2}\big(\pi-\pi^{e},\chi\big)_{D^{k}}+\big(\varphi-\varphi^{e},\psi\big)_{D^{k}}
+\big(\vartheta-\vartheta^{e},\zeta\big)_{D^{k}}
+\omega_{2}\big(\sigma-\sigma^{e},\gamma\big)_{D^{k}}\\
&+ \omega_{1}\big(\pi_{2}-\pi_{2}^{e},\gamma\big)_{D^{k}}+ \lambda_{4}\big(\pi_{2}-\pi_{2}^{e},\gamma\big)_{D^{k}}+\big(\pi_{3}-\pi_{3}^{e},\delta\big)_{D^{k}}
+\big(\pi_{4}-\pi_{4}^{e},\varsigma\big)_{D^{k}}-\omega_{2}\big(\pi-\pi^{e},o\big)_{D^{k}}- \omega_{1}\big(\pi_{1}-\pi_{1}^{e},o\big)_{D^{k}}\\
&- \lambda_{4}\big(\pi_{1}-\pi_{1}^{e},o\big)_{D^{k}}+\big(\pi_{5}-\pi_{5}^{e},\omega\big)_{D^{k}}
+\big(\pi_{6}-\pi_{6}^{e},\kappa\big)_{D^{k}}
-\big(n.(\tau-\tau^{e})^{*},\phi\big)_{\partial D^{k}}-\big(n.(\sigma-\sigma^{e})^{*},\varphi\big)_{\partial D^{k}}\\
&-\big(n.(\vartheta-\vartheta^{e})^{*},\psi\big)_{\partial D^{k}}-\big(n.(\pi-\pi^{e})^{*},\zeta\big)_{\partial D^{k}}-\big(n.(\pi_{4}-\pi_{4}^{e})^{*},\delta\big)_{\partial D^{k}}-\big(n.(\pi_{2}-\pi_{2}^{e})^{*},\varsigma\big)_{\partial D^{k}}\\
&+\big(n.(\pi_{6}-\pi_{6}^{e})^{*},\omega\big)_{\partial D^{k}}
-\big(n.(\pi_{1}-\pi_{1}^{e})^{*},\kappa\big)_{\partial D^{k}}=0.\\
\end{split}
\end{equation}
We take the test functions
\begin{equation}\label{91h}
\begin{split}
&\vartheta_{1}=\pi,\quad\beta_{1}=\epsilon_{1}-\epsilon,\quad \phi=\pi,\quad \varphi=-\vartheta,\quad \chi=\sigma,\quad\beta_{2}=\epsilon_{2}-\varphi,\quad \psi=-\sigma,\quad \zeta=\tau,\\
&\gamma=\pi_{1},\quad\beta_{3}=\epsilon_{3}-\pi_{3},\quad \delta=\pi_{1},\quad \varsigma=-\pi_{6},\quad o=\pi_{2},\quad\beta_{4}=\epsilon_{4}-\pi_{5},\quad \omega=-\pi_{2},\quad \kappa=\pi_{4},\\
\end{split}
\end{equation}
 we obtain
\begin{equation}\label{1bxhhn}
\begin{split}
\big(&\frac{\partial (\pi-\pi^{e})}{\partial t},\pi\big)_{D^{k}}+\big(\frac{\partial (\sigma-\sigma^{e})}{\partial t},\sigma\big)_{D^{k}}+\big(\frac{\partial (\pi_{1}-\pi_{1}^{e})}{\partial t},\pi_{1}\big)_{D^{k}}+\big(\frac{\partial (\pi_{2}-\pi_{2}^{e})}{\partial t},\pi_{2}\big)_{D^{k}}-\big(\Delta_{(\alpha-2)/2}(\epsilon-\epsilon^{e}),\epsilon_{1}-\epsilon\big)_{D^{k}}\\
&-\big(\Delta_{(\alpha-2)/2}(\varphi-\varphi^{e}),\epsilon_{2}-\varphi\big)_{D^{k}}-
\big(\Delta_{(\alpha-2)/2}(\pi_{3}-\pi_{3}^{e}),\epsilon_{3}-\pi_{3}\big)_{D^{k}}
-\big(\Delta_{(\alpha-2)/2}(\pi_{5}-\pi_{5}^{e}),\epsilon_{4}-\pi_{5}\big)_{D^{k}}
\\
&+\lambda_{1}\big(\epsilon_{1}-\epsilon_{1}^{e},\pi\big)_{D^{k}}
-\lambda_{1}\big(\epsilon_{2}-\epsilon_{2}^{e},\sigma\big)_{D^{k}}+
\lambda_{3}\big(\epsilon_{3}-\epsilon_{3}^{e},\pi_{1}\big)_{D^{k}}
-\lambda_{3}\big(\epsilon_{4}-\epsilon_{4}^{e},\pi_{2}\big)_{D^{k}}
+\big(\epsilon_{1}-\epsilon_{1}^{e},\epsilon_{1}-\epsilon\big)_{D^{k}}\\
&+\big(\epsilon_{2}-\epsilon_{2}^{e},\epsilon_{2}-\varphi\big)_{D^{k}}
+\big(\epsilon_{3}-\epsilon_{3}^{e},\epsilon_{3}-\pi_{3}\big)_{D^{k}}
+\big(\epsilon_{4}-\epsilon_{4}^{e},\epsilon_{4}-\pi_{5}\big)_{D^{k}}
+\big(\tau-\tau^{h},\pi_{x}\big)_{D^{k}}-\big(\sigma-\sigma^{h},\vartheta_{x}\big)_{D^{k}}\\
&-\big(\vartheta-\vartheta^{e},\sigma_{x}\big)_{D^{k}}
+\big(\pi-\pi^{e},\tau_{x}\big)_{D^{k}}+\big(\pi_{4}-\pi_{4}^{e},(\pi_{1})_{x}\big)_{D^{k}}
-\big(\pi_{2}-\pi_{2}^{e},(\pi_{6})_{x}\big)_{D^{k}}
+\big(\pi_{1}-\pi_{1}^{e},(\pi_{4})_{x}\big)_{D^{k}}\\
&-\big(\pi_{6}-\pi_{6}^{e},(\pi_{2})_{x}\big)_{D^{k}}+ \omega_{1}\big(\sigma-\sigma^{e},\pi\big)_{D^{k}}+ \omega_{2}\big(\pi_{2}-\pi_{2}^{e},\pi\big)_{D^{k}}+ \lambda_{2}\big(\sigma-\sigma^{e},\pi\big)_{D^{k}}+\big(\epsilon-\epsilon^{e},\pi\big)_{D^{k}}
-\big(\tau-\tau^{e},\vartheta\big)_{D^{k}}
\\
&- \omega_{1}\big(\pi-\pi^{e},\sigma\big)_{D^{k}}-\omega_{2} \big(\pi_{1}-\pi_{1}^{e},\sigma\big)_{D^{k}}- \lambda_{2}\big(\pi-\pi^{e},\sigma\big)_{D^{k}}-\big(\varphi-\varphi^{e},\sigma\big)_{D^{k}}
+\big(\vartheta-\vartheta^{e},\tau\big)_{D^{k}}
+\omega_{2}\big(\sigma-\sigma^{e},\pi_{1}\big)_{D^{k}}\\
&+ \omega_{1}\big(\pi_{2}-\pi_{2}^{e},\pi_{1}\big)_{D^{k}}+ \lambda_{4}\big(\pi_{2}-\pi_{2}^{e},\pi_{1}\big)_{D^{k}}+\big(\pi_{3}-\pi_{3}^{e},\pi_{1}\big)_{D^{k}}
-\big(\pi_{4}-\pi_{4}^{e},\pi_{6}\big)_{D^{k}}-\omega_{2}\big(\pi-\pi^{e},\pi_{2}\big)_{D^{k}}- \omega_{1}\big(\pi_{1}-\pi_{1}^{e},\pi_{2}\big)_{D^{k}}\\
&- \lambda_{4}\big(\pi_{1}-\pi_{1}^{e},\pi_{2}\big)_{D^{k}}-\big(\pi_{5}-\pi_{5}^{e},\pi_{2}\big)_{D^{k}}
+\big(\pi_{6}-\pi_{6}^{e},\pi_{4}\big)_{D^{k}}
-\big(n.(\tau-\tau^{e})^{*},\pi\big)_{\partial D^{k}}+\big(n.(\sigma-\sigma^{e})^{*},\vartheta\big)_{\partial D^{k}}\\
&+\big(n.(\vartheta-\vartheta^{e})^{*},\sigma\big)_{\partial D^{k}}-\big(n.(\pi-\pi^{e})^{*},\tau\big)_{\partial D^{k}}-\big(n.(\pi_{4}-\pi_{4}^{e})^{*},\pi_{1}\big)_{\partial D^{k}}+\big(n.(\pi_{2}-\pi_{2}^{e})^{*},\pi_{6}\big)_{\partial D^{k}}\\
&+\big(n.(\pi_{6}-\pi_{6}^{e})^{*},\pi_{2}\big)_{\partial D^{k}}
-\big(n.(\pi_{1}-\pi_{1}^{e})^{*},\pi_{4}\big)_{\partial D^{k}}=0.\\
\end{split}
\end{equation}
Summing over $k$, simplify  by integration by parts and \eqref{flc}, we get
\begin{equation}\label{tt12hm}
\begin{split}
&\big(\frac{\partial \pi}{\partial t},\pi\big)_{\Omega_{h}}+\big(\frac{\partial \sigma}{\partial t},\sigma\big)_{\Omega_{h}}+\big(\frac{\partial \pi_{1}}{\partial t},\pi_{1}\big)_{D^{k}}+\big(\frac{\partial \pi_{2}}{\partial t},\pi_{2}\big)_{\Omega_{h}}
+\big(\Delta_{(\alpha-2)/2}\epsilon,\epsilon\big)_{\Omega_{h}}+\big(\Delta_{(\alpha-2)/2}\varphi,\varphi\big)_{\Omega_{h}}
\\
&\quad+\big(\Delta_{(\alpha-2)/2}\pi_{3},\pi_{3}\big)_{\Omega_{h}}
+\big(\Delta_{(\alpha-2)/2}\pi_{5},\pi_{5}\big)_{\Omega_{h}}
+\big(\epsilon_{1},\epsilon_{1}\big)_{\Omega_{h}}+\big(\epsilon_{2},\epsilon_{2}\big)_{\Omega_{h}}
+\big(\epsilon_{3},\epsilon_{3}\big)_{\Omega_{h}}
+\big(\epsilon_{4},\epsilon_{4}\big)_{\Omega_{h}}\\
&=\big( (\pi^{e})_{t},\pi\big)_{\Omega_{h}}+\big((\sigma^{e})_{t},\sigma\big)_{\Omega_{h}}+\big( (\pi_{1}^{e})_{t},\pi_{1}\big)_{\Omega_{h}}+\big( (\pi_{2}^{e})_{t},\pi_{2}\big)_{\Omega_{h}}+
\big(\epsilon_{1}^{e},\epsilon_{1}-\epsilon\big)_{\Omega_{h}}
+\big(\epsilon_{2}^{e},\epsilon_{2}-\varphi\big)_{\Omega_{h}}\\
&\quad+\big(\epsilon_{3}^{e},\epsilon_{3}-\pi_{3}\big)_{\Omega_{h}}
+\big(\epsilon_{4}^{e},\epsilon_{4}-\pi_{5}\big)_{\Omega_{h}}
-\big(\Delta_{(\alpha-2)/2}\epsilon^{e},\epsilon_{1}-\epsilon\big)_{\Omega_{h}}
-\big(\Delta_{(\alpha-2)/2}\varphi^{e},\epsilon_{2}-\varphi\big)_{\Omega_{h}}\\
&\quad-
\big(\Delta_{(\alpha-2)/2}\pi_{3}^{e},\epsilon_{3}-\pi_{3}\big)_{\Omega_{h}}
-\big(\Delta_{(\alpha-2)/2}\pi_{5}^{e},\epsilon_{4}-\pi_{5}\big)_{\Omega_{h}}
+\omega_{1}\big(\sigma^{e},\pi\big)_{\Omega_{h}}- \big(\pi^{e},\sigma\big)_{\Omega_{h}}+\big(\epsilon^{e},\pi\big)_{\Omega_{h}}
-\big(\varphi^{e},\sigma\big)_{\Omega_{h}} \\
&\quad+ \lambda_{2}\big(\sigma^{e},\pi\big)_{\Omega_{h}}+ \omega_{2}\big(\pi_{2}^{e},\pi\big)_{\Omega_{h}}+ \big(\sigma^{e},\pi\big)_{\Omega_{h}}+\big(\epsilon^{e},\pi\big)_{\Omega_{h}}
-\big(\tau^{e},\vartheta\big)_{\Omega_{h}}
+\lambda_{1}\big(\epsilon_{1}^{e},\pi\big)_{\Omega_{h}}-\lambda_{1}\big(\epsilon_{2}^{e},\sigma\big)_{\Omega_{h}}\\
&\quad- \omega_{2}\big(\pi^{e},\sigma\big)_{\Omega_{h}}- \omega_{1}\big(\pi_{1}^{e},\sigma\big)_{\Omega_{h}}- \lambda_{2}\big(\pi^{e},\sigma\big)_{\Omega_{h}}-\big(\varphi^{e},\sigma\big)_{\Omega_{h}}
+\big(\vartheta^{e},\tau\big)_{\Omega_{h}}
+\omega_{2}\big(\sigma^{e},\pi_{1}\big)_{\Omega_{h}}
+\omega_{1}\big(\pi_{2}^{e},\pi_{1}\big)_{\Omega_{h}}\\
&\quad
-\big(\pi_{4}^{e},\pi_{6}\big)_{\Omega_{h}}-\omega_{2}\big(\pi^{e},\pi_{2}\big)_{\Omega_{h}}- \omega_{1}\big(\pi_{1}^{e},\pi_{2}\big)_{\Omega_{h}}- \lambda_{2}\big(\pi_{1}^{e},\pi_{2}\big)_{\Omega_{h}}-\big(\pi_{5}^{e},\pi_{2}\big)_{\Omega_{h}}
+\big(\pi_{6}^{e},\pi_{4}\big)_{\Omega_{h}}\\
&\quad+\lambda_{3}\big(\epsilon_{3}^{e},\pi_{1}\big)_{\Omega_{h}}-\lambda_{3}\big(\epsilon_{4}^{e},\pi_{2}\big)_{\Omega_{h}}
+ \lambda_{4}\big(\pi_{2}^{e},\pi_{1}\big)_{\Omega_{h}}+\big(\pi_{3}^{e},\pi_{1}\big)_{\Omega_{h}}\\
&\quad+\big(\Delta_{(\alpha-2)/2}\epsilon,\epsilon_{1}\big)_{\Omega_{h}}
+\big(\Delta_{(\alpha-2)/2}\varphi,\epsilon_{2}\big)_{\Omega_{h}}+
\big(\Delta_{(\alpha-2)/2}\pi_{3},\epsilon_{3}\big)_{\Omega_{h}}
+\big(\Delta_{(\alpha-2)/2}\pi_{5},\epsilon_{4}\big)_{\Omega_{h}}
+\big(\epsilon_{1},\epsilon\big)_{\Omega_{h}}\\
&\quad+\big(\epsilon_{2},\varphi\big)_{\Omega_{h}}
+\big(\epsilon_{3},\pi_{3}\big)_{\Omega_{h}}
+\big(\epsilon_{4},\pi_{5}\big)_{\Omega_{h}}
-\big(\epsilon,\pi\big)_{\Omega_{h}}-\big(\pi_{3},\pi_{1}\big)_{\Omega_{h}}
+\big(\pi_{5},\pi_{2}\big)_{\Omega_{h}}+\big(\varphi,\sigma\big)_{\Omega_{h}}\\
&\quad-\lambda_{1}\big(\epsilon_{1},\pi\big)_{\Omega_{h}}+\lambda_{1}\big(\epsilon_{2},\sigma\big)_{\Omega_{h}}-
\lambda_{3}\big(\epsilon_{3},\pi_{1}\big)_{\Omega_{h}}+\lambda_{3}\big(\epsilon_{4},\pi_{2}\big)_{\Omega_{h}}\\
&\quad-\sum_{k=1}^{K}((\tau^{e})^{+}[\pi])_{k+\frac{1}{2}}
+\sum_{k=1}^{K}((\sigma^{e})^{-}[\vartheta])_{k+\frac{1}{2}}+
\sum_{k=1}^{K}((\vartheta^{e})^{+}[\sigma])_{k+\frac{1}{2}}-\sum_{k=1}^{K}((\pi^{e})^{-}[\tau])_{k+\frac{1}{2}}\\
&\quad-\sum_{k=1}^{K}((\pi_{4}^{e})^{+}[\pi_{1}])_{k+\frac{1}{2}}
+\sum_{k=1}^{K}((\pi_{2}^{e})^{-}[\pi_{6}])_{k+\frac{1}{2}}+
\sum_{k=1}^{K}((\pi_{6}^{e})^{+}[\pi_{2}])_{k+\frac{1}{2}}
-\sum_{k=1}^{K}((\pi_{1}^{e})^{-}[\pi_{4}])_{k+\frac{1}{2}}\\
&\quad+\big(\tau^{e},\pi_{x}\big)_{\Omega_{h}}-\big(\sigma^{e},\vartheta_{x}\big)_{\Omega_{h}}
-\big(\vartheta^{e},\sigma_{x}\big)_{D^{k}}
+\big(\pi^{e},\tau_{x}\big)_{\Omega_{h}}+\big(\pi_{4}^{e},(\pi_{1})_{x}\big)_{\Omega_{h}}
-\big(\pi_{2}^{e},(\pi_{6})_{x}\big)_{\Omega_{h}}
+\big(\pi_{1}^{e},(\pi_{4})_{x}\big)_{\Omega_{h}}\\
&\quad-\big(\pi_{6}^{e},(\pi_{2})_{x}\big)_{\Omega_{h}}=T_{1}+T_{2}+T_{3}+T_{4}.
\end{split}
\end{equation}
Now, we estimate $T_{i}$ term by term.
\begin{equation}\label{tt12h}
\begin{split}
T_{1}=&\big( (\pi^{e})_{t},\pi\big)_{\Omega_{h}}+\big((\sigma^{e})_{t},\sigma\big)_{\Omega_{h}}+\big( (\pi_{1}^{e})_{t},\pi_{1}\big)_{\Omega_{h}}+\big( (\pi_{2}^{e})_{t},\pi_{2}\big)_{\Omega_{h}}+
\big(\epsilon_{1}^{e},\epsilon_{1}-\epsilon\big)_{\Omega_{h}}
+\big(\epsilon_{2}^{e},\epsilon_{2}-\varphi\big)_{\Omega_{h}}
\\
&+\big(\epsilon_{3}^{e},\epsilon_{3}-\pi_{3}\big)_{\Omega_{h}}
+\big(\epsilon_{4}^{e},\epsilon_{4}-\pi_{5}\big)_{\Omega_{h}}
-\big(\Delta_{(\alpha-2)/2}\epsilon^{e},\epsilon_{1}-\epsilon\big)_{\Omega_{h}}
-\big(\Delta_{(\alpha-2)/2}\varphi^{e},\epsilon_{2}-\varphi\big)_{\Omega_{h}}\\
&-
\big(\Delta_{(\alpha-2)/2}\pi_{3}^{e},\epsilon_{3}-\pi_{3}\big)_{\Omega_{h}}
-\big(\Delta_{(\alpha-2)/2}\pi_{5}^{e},\epsilon_{4}-\pi_{5}\big)_{\Omega_{h}}
+\big(\sigma^{e},\pi\big)_{\Omega_{h}}- \big(\pi^{e},\sigma\big)_{\Omega_{h}}+\big(\epsilon^{e},\pi\big)_{\Omega_{h}}
-\big(\varphi^{e},\sigma\big)_{\Omega_{h}} \\
&+ \lambda_{2}\big(\sigma^{e},\pi\big)_{\Omega_{h}}+ \omega_{2}\big(\pi_{2}^{e},\pi\big)_{\Omega_{h}}+ \omega_{1}\big(\sigma^{e},\pi\big)_{\Omega_{h}}+\big(\epsilon^{e},\pi\big)_{\Omega_{h}}
+\lambda_{1}\big(\epsilon_{1}^{e},\pi\big)_{D^{k}}-\lambda_{1}\big(\epsilon_{2}^{e},\sigma\big)_{\Omega_{h}}\\
&- \omega_{2}\big(\pi^{e},\sigma\big)_{\Omega_{h}}- \omega_{1}\big(\pi_{1}^{e},\sigma\big)_{\Omega_{h}}- \lambda_{2}\big(\pi^{e},\sigma\big)_{\Omega_{h}}-\big(\varphi^{e},\sigma\big)_{\Omega_{h}}
+\omega_{2}\big(\sigma^{e},\pi_{1}\big)_{\Omega_{h}}
+\omega_{1} \big(\pi_{2}^{e},\pi_{1}\big)_{\Omega_{h}}\\
&
-\omega_{2}\big(\pi^{e},\pi_{2}\big)_{\Omega_{h}}- \omega_{1}\big(\pi_{1}^{e},\pi_{2}\big)_{\Omega_{h}}- \lambda_{2}\big(\pi_{1}^{e},\pi_{2}\big)_{\Omega_{h}}-\big(\pi_{5}^{e},\pi_{2}\big)_{\Omega_{h}}
+\lambda_{3}\big(\epsilon_{3}^{e},\pi_{1}\big)_{\Omega_{h}}\\
&-\lambda_{3}\big(\epsilon_{4}^{e},\pi_{2}\big)_{\Omega_{h}}
+ \lambda_{4}\big(\pi_{2}^{e},\pi_{1}\big)_{\Omega_{h}}+\big(\pi_{3}^{e},\pi_{1}\big)_{\Omega_{h}}.
\\
\end{split}
\end{equation}
Employing Young's inequality, we obtain
\begin{equation}\label{tt12h1}
\begin{split}
T_{1}\leq &c_{12}\|\epsilon\|^{2}_{L^{2}(\Omega_{h})}
+c_{11}\|\pi\|^{2}_{L^{2}(\Omega_{h})}+c_{10}\|\sigma\|^{2}_{L^{2}(\Omega_{h})}+
c_{9}\|\pi_{1}\|^{2}_{L^{2}(\Omega_{h})}+c_{8}\|\pi_{2}\|^{2}_{L^{2}(\Omega_{h})}+
c_{7}\|\pi_{3}\|^{2}_{L^{2}(\Omega_{h})}\\
&+c_{6}\|\pi_{5}\|^{2}_{L^{2}(\Omega_{h})}+
c_{1}\|\epsilon_{1}\|^{2}_{L^{2}(\Omega_{h})}+c_{2}\|\epsilon_{2}\|^{2}_{L^{2}(\Omega_{h})}
+c_{3}\|\epsilon_{3}\|^{2}_{L^{2}(\Omega_{h})}\\
&+c_{4}\|\epsilon_{4}\|^{2}_{L^{2}(\Omega_{h})}
+c_{5}\|\varphi\|^{2}_{L^{2}(\Omega_{h})}+Ch^{2N+2},\\
\end{split}
\end{equation}
and
\begin{equation}\label{tt12h}
\begin{split}
T_{2}=&\big(\Delta_{(\alpha-2)/2}\epsilon,\epsilon_{1}\big)_{\Omega_{h}}
+\big(\Delta_{(\alpha-2)/2}\varphi,\epsilon_{2}\big)_{\Omega_{h}}+
\big(\Delta_{(\alpha-2)/2}\pi_{3},\epsilon_{3}\big)_{\Omega_{h}}
+\big(\Delta_{(\alpha-2)/2}\pi_{5},\epsilon_{4}\big)_{\Omega_{h}}
+\big(\epsilon_{1},\epsilon\big)_{\Omega_{h}}\\
&+\big(\epsilon_{2},\varphi\big)_{\Omega_{h}}
+\big(\epsilon_{3},\pi_{3}\big)_{\Omega_{h}}
+\big(\epsilon_{4},\pi_{5}\big)_{\Omega_{h}}
-\big(\epsilon,\pi\big)_{\Omega_{h}}-\big(\pi_{3},\pi_{1}\big)_{\Omega_{h}}
+\big(\pi_{5},\pi_{2}\big)_{\Omega_{h}}+\big(\varphi,\sigma\big)_{\Omega_{h}}\\
&-\lambda_{1}\big(\epsilon_{1},\pi\big)_{\Omega_{h}}+\lambda_{1}\big(\epsilon_{2},\sigma\big)_{\Omega_{h}}-
\lambda_{3}\big(\epsilon_{3},\pi_{1}\big)_{\Omega_{h}}+\lambda_{3}\big(\epsilon_{4},\pi_{2}\big)_{\Omega_{h}}.\\
\end{split}
\end{equation}
Employing Young's inequality and  Lemma \ref{lga2}, we obtain
\begin{equation}\label{tt12h2}
\begin{split}
T_{2}\leq &c_{12}\|\epsilon\|^{2}_{L^{2}(\Omega_{h})}
+c_{11}\|\pi\|^{2}_{L^{2}(\Omega_{h})}+c_{10}\|\sigma\|^{2}_{L^{2}(\Omega_{h})}+
c_{9}\|\pi_{1}\|^{2}_{L^{2}(\Omega_{h})}+c_{8}\|\pi_{2}\|^{2}_{L^{2}(\Omega_{h})}+
c_{7}\|\pi_{3}\|^{2}_{L^{2}(\Omega_{h})}\\
&+c_{6}\|\pi_{5}\|^{2}_{L^{2}(\Omega_{h})}+
c_{1}\|\epsilon_{1}\|^{2}_{L^{2}(\Omega_{h})}+c_{2}\|\epsilon_{2}\|^{2}_{L^{2}(\Omega_{h})}
+c_{3}\|\epsilon_{3}\|^{2}_{L^{2}(\Omega_{h})}+c_{4}\|\epsilon_{4}\|^{2}_{L^{2}(\Omega_{h})}
+c_{5}\|\varphi\|^{2}_{L^{2}(\Omega_{h})}.\\
\end{split}
\end{equation}
and
\begin{equation}\label{tt12h}
\begin{split}
T_{3}=&-\sum_{k=1}^{K}((\tau^{e})^{+}[\pi])_{k+\frac{1}{2}}
+\sum_{k=1}^{K}((\sigma^{e})^{-}[\vartheta])_{k+\frac{1}{2}}+
\sum_{k=1}^{K}((\vartheta^{e})^{+}[\sigma])_{k+\frac{1}{2}}-\sum_{k=1}^{K}((\pi^{e})^{-}[\tau])_{k+\frac{1}{2}}\\
&-\sum_{k=1}^{K}((\pi_{4}^{e})^{+}[\pi_{1}])_{k+\frac{1}{2}}
+\sum_{k=1}^{K}((\pi_{2}^{e})^{-}[\pi_{6}])_{k+\frac{1}{2}}+
\sum_{k=1}^{K}((\pi_{6}^{e})^{+}[\pi_{2}])_{k+\frac{1}{2}}-\sum_{k=1}^{K}((\pi_{1}^{e})^{-}[\pi_{4}])_{k+\frac{1}{2}}.\\
\end{split}
\end{equation}
and
\begin{equation}\label{tt12h}
\begin{split}
T_{4}=&\big(\tau^{h},\pi_{x}\big)_{D^{k}}-\big(\sigma^{h},\vartheta_{x}\big)_{\Omega_{h}}
-\big(\vartheta^{e},\sigma_{x}\big)_{\Omega_{h}}
+\big(\pi^{e},\tau_{x}\big)_{\Omega_{h}}+\big(\pi_{4}^{e},(\pi_{1})_{x}\big)_{\Omega_{h}}
-\big(\pi_{2}^{e},(\pi_{6})_{x}\big)_{\Omega_{h}}\\
&+\big(\pi_{1}^{e},(\pi_{4})_{x}\big)_{\Omega_{h}}-\big(\pi_{6}^{e},(\pi_{2})_{x}\big)_{\Omega_{h}}
+\big(\pi_{6}^{e},\pi_{4}\big)_{\Omega_{h}}-\big(\pi_{4}^{e},\pi_{6}\big)_{\Omega_{h}}
-\big(\tau^{e},\vartheta\big)_{\Omega_{h}}+\big(\vartheta^{e},\tau\big)_{\Omega_{h}}.\\
\end{split}
\end{equation}
Using the definition of the numerical traces, \eqref{flc}, and the definitions of the projections $\mathcal{P}^{+},\mathcal{P}^{-}$ \eqref{prh}, we get
\begin{equation}\label{tt12h3}
\begin{split}
T_{3}=T_{4}=0.
\end{split}
\end{equation}
Combining \eqref{tt12h1}, \eqref{tt12h2},  \eqref{tt12h3} and \eqref{tt12hm}, we obtain
\begin{equation}\label{tt12hmm}
\begin{split}
&\big(\frac{\partial \pi}{\partial t},\pi\big)_{\Omega_{h}}+\big(\frac{\partial \sigma}{\partial t},\sigma\big)_{\Omega_{h}}+\big(\frac{\partial \pi_{1}}{\partial t},\pi_{1}\big)_{D^{k}}+\big(\frac{\partial \pi_{2}}{\partial t},\pi_{2}\big)_{\Omega_{h}}
+\big(\Delta_{(\alpha-2)/2}\epsilon,\epsilon\big)_{\Omega_{h}}+\big(\Delta_{(\alpha-2)/2}\varphi,\varphi\big)_{\Omega_{h}}
\\
&\quad+ \big(\Delta_{(\alpha-2)/2}\pi_{3},\pi_{3}\big)_{\Omega_{h}}
+\big(\Delta_{(\alpha-2)/2}\pi_{5},\pi_{5}\big)_{\Omega_{h}}
+\big(\epsilon_{1},\epsilon_{1}\big)_{\Omega_{h}}+\big(\epsilon_{2},\epsilon_{2}\big)_{\Omega_{h}}
+\big(\epsilon_{3},\epsilon_{3}\big)_{\Omega_{h}}
+\big(\epsilon_{4},\epsilon_{4}\big)_{\Omega_{h}}\\
&\leq c_{9}\|\epsilon\|^{2}_{L^{2}(\Omega_{h})}
+c_{5}\|\pi\|^{2}_{L^{2}(\Omega_{h})}+c_{6}\|\sigma\|^{2}_{L^{2}(\Omega_{h})}+
c_{7}\|\pi_{1}\|^{2}_{L^{2}(\Omega_{h})}+c_{8}\|\pi_{2}\|^{2}_{L^{2}(\Omega_{h})}+
c_{12}\|\pi_{3}\|^{2}_{L^{2}(\Omega_{h})}\\
&\quad+c_{11}\|\pi_{5}\|^{2}_{L^{2}(\Omega_{h})}+
c_{1}\|\epsilon_{1}\|^{2}_{L^{2}(\Omega_{h})}+c_{2}\|\epsilon_{2}\|^{2}_{L^{2}(\Omega_{h})}
+c_{3}\|\epsilon_{3}\|^{2}_{L^{2}(\Omega_{h})}+c_{4}\|\epsilon_{4}\|^{2}_{L^{2}(\Omega_{h})}
+c_{10}\|\varphi\|^{2}_{L^{2}(\Omega_{h})}+Ch^{2N+2}.\\
\end{split}
\end{equation}
Recalling Lemma \ref{lg}, we get

\begin{equation}\label{tt12hmb}
\begin{split}
&\big(\frac{\partial \pi}{\partial t},\pi\big)_{\Omega_{h}}+\big(\frac{\partial \sigma}{\partial t},\sigma\big)_{\Omega_{h}}+\big(\frac{\partial \pi_{1}}{\partial t},\pi_{1}\big)_{D^{k}}+\big(\frac{\partial \pi_{2}}{\partial t},\pi_{2}\big)_{\Omega_{h}}
+\big(\epsilon_{1},\epsilon_{1}\big)_{\Omega_{h}}+\big(\epsilon_{2},\epsilon_{2}\big)_{\Omega_{h}}
+\big(\epsilon_{3},\epsilon_{3}\big)_{\Omega_{h}}
+\big(\epsilon_{4},\epsilon_{4}\big)_{\Omega_{h}}\\
&\leq c_{5}\|\pi\|^{2}_{L^{2}(\Omega_{h})}+c_{6}\|\sigma\|^{2}_{L^{2}(\Omega_{h})}+
c_{7}\|\pi_{1}\|^{2}_{L^{2}(\Omega_{h})}+c_{8}\|\pi_{2}\|^{2}_{L^{2}(\Omega_{h})}
\\
&\,\,\,\,+
c_{1}\|\epsilon_{1}\|^{2}_{L^{2}(\Omega_{h})}+c_{2}\|\epsilon_{2}\|^{2}_{L^{2}(\Omega_{h})}
+c_{3}\|\epsilon_{3}\|^{2}_{L^{2}(\Omega_{h})}+c_{4}\|\epsilon_{4}\|^{2}_{L^{2}(\Omega_{h})}
+Ch^{2N+2},\\
\end{split}
\end{equation}
provided $c_{i},\,\,i=1,2,...,8$ are sufficiently small such that $c_{i}\leq1$, we obtain
\begin{equation}\label{tt12hmv}
\begin{split}
&\big(\frac{\partial \pi}{\partial t},\pi\big)_{\Omega_{h}}+\big(\frac{\partial \sigma}{\partial t},\sigma\big)_{\Omega_{h}}+\big(\frac{\partial \pi_{1}}{\partial t},\pi_{1}\big)_{D^{k}}+\big(\frac{\partial \pi_{2}}{\partial t},\pi_{2}\big)_{\Omega_{h}}
\\
&\leq \|\pi\|^{2}_{L^{2}(\Omega_{h})}+\|\sigma\|^{2}_{L^{2}(\Omega_{h})}+
\|\pi_{1}\|^{2}_{L^{2}(\Omega_{h})}+\|\pi_{2}\|^{2}_{L^{2}(\Omega_{h})}
+Ch^{2N+2}.\\
\end{split}
\end{equation}
An integration in $t$ plus the standard approximation theory then
gives the desired error estimates.\\
\section{Numerical examples}\label{sc5} In this section we will present several numerical examples to illustrate the previous theoretical results. Before that, we adopt the nodal discontinuous Galerkin methods for the full spatial discretization   using   a high-order nodal basis set of orthonormal Lagrange-Legendre polynomials of arbitrary order in space on each element of   computational  domain as a more suitable and computationally stable approach As shown by  Aboelenen and El-Hawary \cite{cann}. We use the high-order Runge-Kutta time discretizations \cite{Cockburn1999}, when the polynomials are of
degree $N$, a higher-order accurate Runge-Kutta (RK) method must be used in order to guarantee
that the scheme is stable. In this paper we use a fourth-order non-Total variation diminishing (TVD) Runge-Kutta scheme \cite{Gottlieb:1998:TVD:279724.279737}. Numerical experiments demonstrate its numerical stability
\begin{equation}\label{a1}
\begin{split}
\frac{\partial \mathbf{u}_{h}}{\partial t}=\mathcal{F}(\mathbf{u}_{h},t),
\end{split}
\end{equation}
where $\mathbf{u}_{h}$ is the vector of unknowns, we can use the standard fourth-order
four stage explicit RK method (ERK)
\begin{equation}\label{a1}
\begin{split}
&\mathbf{k}^{1}=\mathcal{F}(\mathbf{u}_{h}^{n},t^{n}),\\
&\mathbf{k}^{2}=\mathcal{F}(\mathbf{u}_{h}^{n}+\frac{1}{2}\Delta t\mathbf{k}^{1},t^{n}+\frac{1}{2}\Delta t),\\
&\mathbf{k}^{3}=\mathcal{F}(\mathbf{u}_{h}^{n}+\frac{1}{2}\Delta t\mathbf{k}^{2},t^{n}+\frac{1}{2}\Delta t),\\
&\mathbf{k}^{4}=\mathcal{F}(\mathbf{u}_{h}^{n}+\Delta t\mathbf{k}^{3},t^{n}+\Delta t),\\
&\mathbf{u}_{h}^{n+1}=\mathbf{u}_{h}^{n}+\frac{1}{6}(\mathbf{k}^{1}+2\mathbf{k}^{2}+2\mathbf{k}^{3}+\mathbf{k}^{4}),
\end{split}
\end{equation}
to advance from $\mathbf{u}_{h}^{n}$ to $\mathbf{u}_{h}^{n+1}$, separated by the time step, $\Delta t$. In our examples, the condition $\Delta t\leq C \Delta x^{\alpha}_{min}\,\,\, (0<C<1)$ is used to ensure stability.
\begin{exmp}\label{ex1} As the first example, we consider the linear fractional Schr\"{o}dinger equation
\begin{equation}\label{sch1}
\begin{split}
&i\frac{\partial u}{\partial t}- \lambda_{1}(-\Delta)^{\frac{\alpha}{2}}u+u=g(x,t),\quad x\in[0,1],\quad t\in(0,0.5],\\
&u(x,0) = u_{0}(x),
\end{split}
\end{equation}

with the initial condition $u_{0}(x)=x^{6}$ and the corresponding forcing term $g(x,t)$  is of the form
\begin{equation}\label{91}
\begin{split}
g(x,t)=e^{-it}\bigg(iu_{0}(x)-\lambda(-\Delta)^{\frac{\alpha}{2}}u_{0}(x)+u_{0}(x)\bigg),
\end{split}
\end{equation}
\end{exmp}
to obtain an exact solution $u(x,t)=e^{-it}x^{6}$
with $\nu=1.2,\,\lambda=\frac{\Gamma(8-\nu)}{2\Gamma(8)}$. The errors
and order of convergence are listed in Table \ref{Tab:a}, confirming optimal $O(h^{N+1})$ order of
convergence across.

\begin{table}[!htb]
    \centering
\begin{center}
 \begin{tabular}{|c|| c c||c|| c c ||c||c c |}
  \hline
 \hline

 \hline
  N&\multicolumn{8}{|c|}{   N=1\quad \quad \quad\quad\quad \quad\qquad \qquad N=2\quad\quad\quad \qquad\quad \quad  \qquad \quad N=3}\\
\hline\hline
  K & $L^{2}$-Error & order & K &$L^{2}$-Error & order&K&$L^{2}$-Error& order \\ [0.5ex]

 \hline
64 & 1.57e-02 & -&       35&8.47e-05&-&         20&1.59e-05&-  \\

74 & 1.24e-02 &1.63&   45&3.97e-05&3.0&    40&9.82e-07&4.02\\
 84&9.2e-03&2.33&     90&5.67e-06&2.81&    60&2.14e-07&3.75\\

 \hline
 \hline

\end{tabular}
\end{center}
\caption{ $L^{2}$-Error and order of convergence for Example \ref{ex1} with $K$ elements and polynomial order $N$.}\label{Tab:a}
\end{table}

\begin{exmp}\label{ex2} Consider  the following nonlinear fractional Schr\"{o}dinger equation
\begin{equation}\label{sch1}
\begin{split}
&i\frac{\partial u}{\partial t}- \lambda(-\Delta)^{\frac{\alpha}{2}}u+|u|^{2}u=g(x,t),\quad x\in[0,1],\quad t\in(0,0.5],\\
&u(x,0) = u_{0}(x),
\end{split}
\end{equation}

with the initial condition $u_{0}(x)=x^{7}$ and the corresponding forcing term $g(x,t)$  is of the form
\begin{equation}\label{91}
\begin{split}
g(x,t)=e^{-it}\bigg(iu_{0}(x)-\lambda(-\Delta)^{\frac{\alpha}{2}}u_{0}(x)+(u_{0}(x))^{3}\bigg).
\end{split}
\end{equation}
\end{exmp}
The exact solution $u(x,t)=e^{-it}x^{7}$ with  $\nu=1.1,\,\lambda=\frac{\Gamma(8-\nu)}{\Gamma(8)}$. The errors and order of convergence are listed in Table \ref{Tab:b}, confirming optimal $O(h^{N+1})$ order of
convergence across.

\begin{table}[!htb]
    \centering
\begin{center}
 \begin{tabular}{|c|| c c||c|| c c ||c||c c |}
  \hline
 \hline

 \hline
  N&\multicolumn{8}{|c|}{   N=1\quad \quad \quad\quad\quad \quad\qquad \qquad N=2\quad\quad\quad \qquad\quad \quad  \qquad \quad N=3}\\
\hline\hline
  K & $L^{2}$-Error & order & K &$L^{2}$-Error & order&K&$L^{2}$-Error& order \\ [0.5ex]

 \hline
120 & 1.41e-01 & -&       60&1.52e-04&-&         40&7.02e-06&-  \\

135 & 1.09e-02 &2.15&   80&6.54e-05&2.89 &    70&7.62e-07&3.97\\
 150&8.9e-03&1.92&     120&1.78e-05&3.22&    90&2.6e-07&4.28\\

 \hline
 \hline

\end{tabular}
\end{center}
\caption{ $L^{2}$-Error and order of convergence for Example \ref{ex2} with $K$ elements and polynomial order $N$.}\label{Tab:b}
\end{table}

\begin{figure}[h!]
\begin{center}
  \vspace{4mm}
\begin{tabular}{cc}
\hspace{1.5cm}
\begin{overpic}[width=5in]{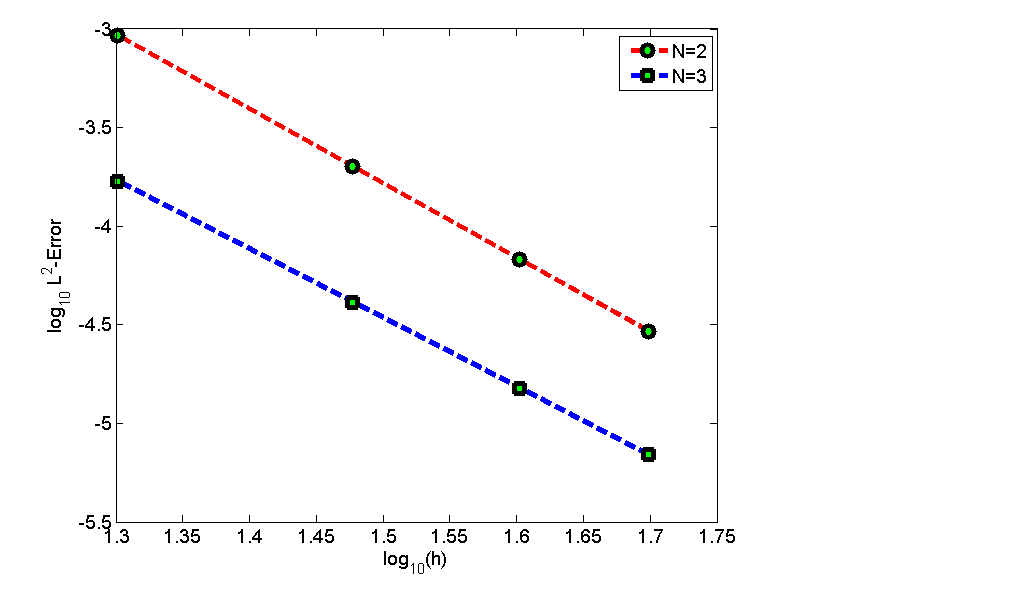}
\end{overpic}

\end{tabular}

\end{center}
\caption[]{\small{ Convergence tests of  \eqref{ex4} with different values of $N$ and $K$.
}} \label{fig:canh}
\end{figure}

\begin{figure}
\begin{center}
  \vspace{5mm}
\begin{tabular}{cc}
\hspace{-0.5cm}
\begin{overpic}[width=4.0in]{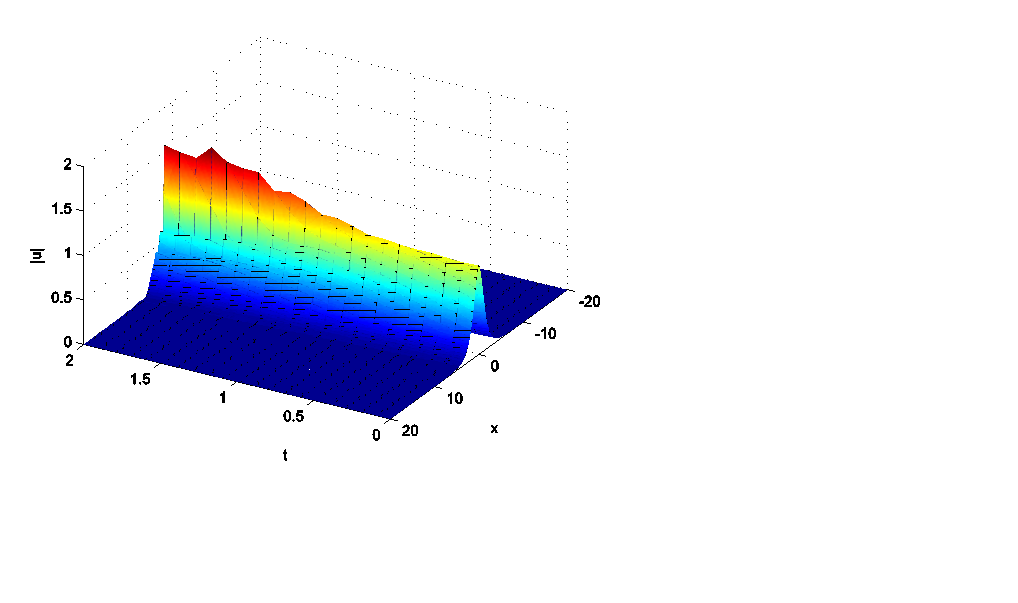}
\put(30,10.75) {\scriptsize \large{$\alpha=1.1$}}

\end{overpic}
 &\hspace{-4.0cm}
\begin{overpic}[width=4.0in]{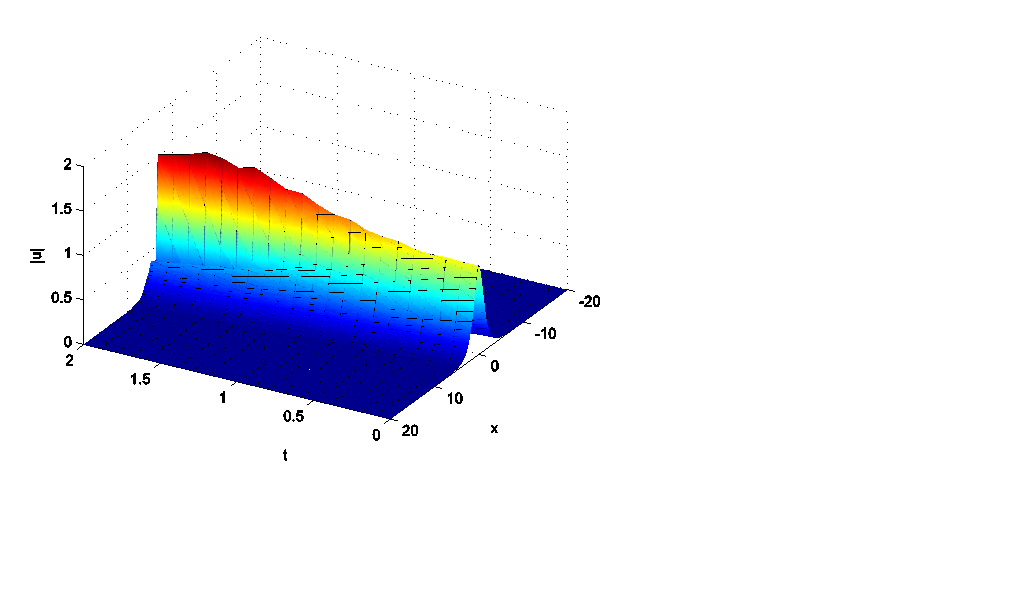}
\put(30,10.75) {\scriptsize \large{$\alpha=1.4$}}

\end{overpic}
%
\\
\\
\\
\hspace{-0.5cm}
\begin{overpic}[width=4in]{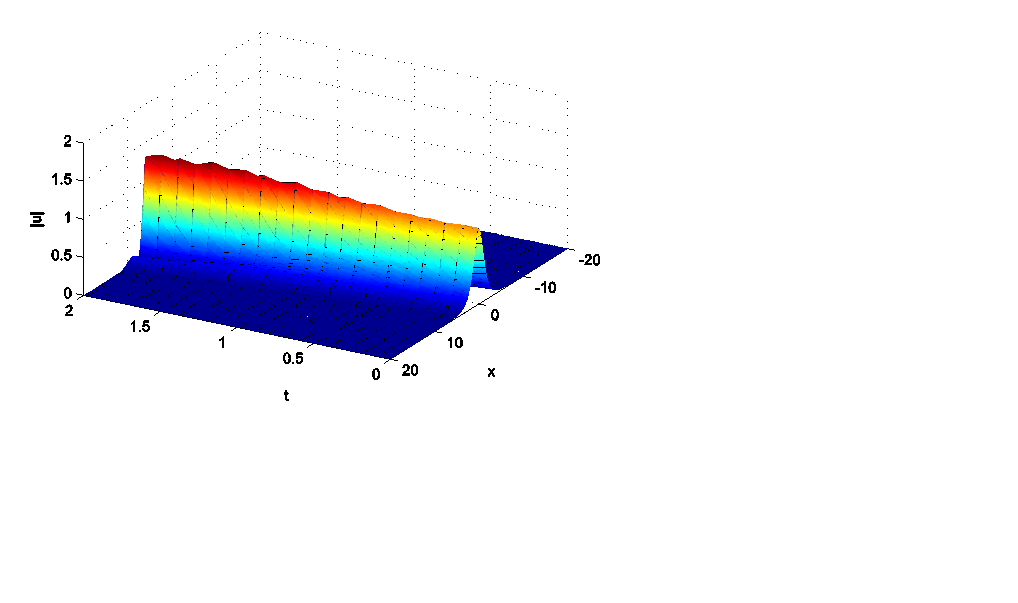}
\put(30,10.75) {\scriptsize \large{$\alpha=1.8$}}
\end{overpic}
 &\hspace{-4.0cm}
\begin{overpic}[width=4in]{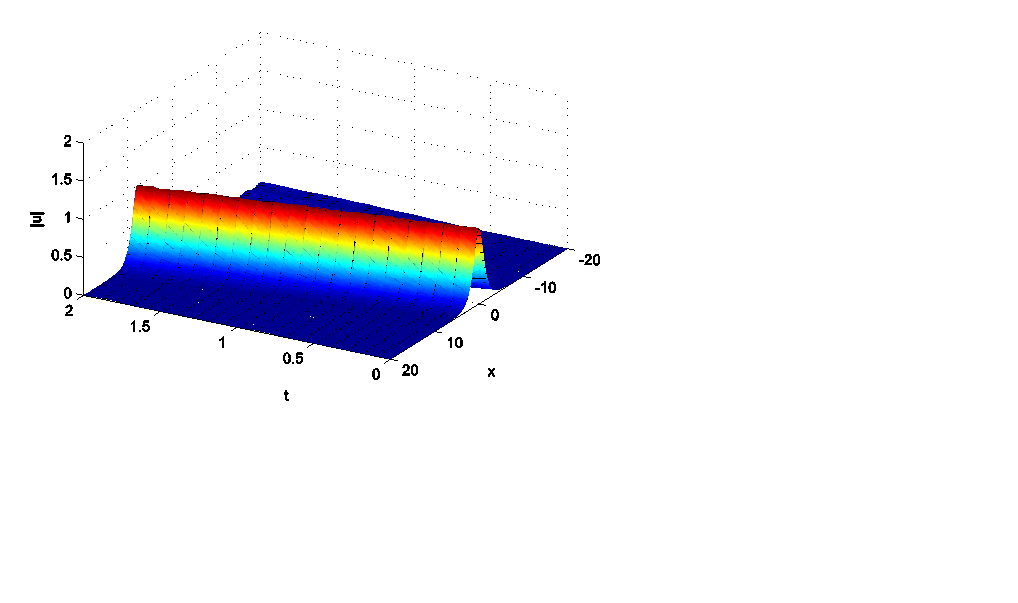}
\put(30,10.75) {\scriptsize \large{$\alpha=2$}}

\end{overpic}

\end{tabular}
\caption[]{{ Numerical results for the nonlinear fractional Schr\"{o}dinger equation   in Example \ref{ex4n}.}}\label{fig:fig2}
\end{center}
\end{figure}

\begin{exmp}  \label{ex4} We consider the nonlinear fractional Schr\"{o}dinger equation
\begin{equation}\label{sch1}
\begin{split}
&i\frac{\partial u}{\partial t}- \lambda(-\Delta)^{\frac{\alpha}{2}}u+|u|^{2}u=g(x,t),\quad x\in[-1,1],\quad t\in(0,0.5],\\
&u(x,0) = u_{0}(x),
\end{split}
\end{equation}

with the initial condition $u_{0}(x)=(x^{2}-1)^{6}$ and the corresponding forcing term $g(x,t)$  is of the form
\begin{equation}\label{91}
\begin{split}
g(x,t)=e^{-it}\bigg(iu_{0}(x)-\lambda(-\Delta)^{\frac{\alpha}{2}}u_{0}(x)+(u_{0}(x))^{3}\bigg),
\end{split}
\end{equation}
\end{exmp}
to obtain an exact solution $u(x,t)=e^{-it}(x^{2}-1)^{6}$ with  $\nu=1.5,\,\lambda=\frac{0.2\Gamma(13-\nu)}{\Gamma(13)}$. We consider cases with $N = 2, 3$ and $K = 20, 30, 40,50$. The numerical orders of convergence
are shown in Figure~\ref{fig:canh}, showing an
$O(h^{N+1})$ convergence rate for all orders.
\begin{exmp}  \label{ex4n} We consider the nonlinear fractional Schr\"{o}dinger equation \eqref{sch1vn} with
 initial condition,
\begin{equation}\label{91}
\begin{split}
u(x,0)=e^{2ix}sech(x),
\end{split}
\end{equation}
\end{exmp}
 with parameters   $\lambda_{1}=\lambda_{2}=1$ and $\,x\in[-20,20]$. We consider cases with $N =  2$ and $K = 80$ and  solve
the equation for several different values of $\alpha$. The numerical solution $u_{h}(x, t)$ for $\alpha=1.1,\,1.4,\,1.8,\,2.0$ is shown in Figure \ref{fig:fig2}.  We observe that
the order $\alpha$ will affect the shape of the soliton case.  When $\alpha$ becomes smaller, the shape of the soliton will change more quickly. This property of the fractional Schr\"{o}dinger equation can be used in physics to modify the shape of wave without change of the nonlinearity and dispersion effects. The numerical solutions of the fractional equation are convergent to the solutions of the classical non-fractional equation when $\alpha$ tends to $2$.


\begin{exmp}\label{ex5} Consider  the  linear coupled fractional Schr\"{o}dinger equations
\begin{equation}\label{sch1}
\begin{split}
&i\frac{\partial u_{1}(x,t)}{\partial t}- \lambda_{1}(-\Delta)^{\frac{\alpha}{2}}u_{1}(x,t)+u_{2}(x,t)+2u_{1}(x,t)=g_{1}(x,t),\,\, x\in[0,1],\,\, t\in(0,0.5],\\
&i\frac{\partial u_{2}(x,t)}{\partial t}- \lambda_{2}(-\Delta)^{\frac{\alpha}{2}}u_{2}(x,t)+2u_{2}(x,t)-u_{1}(x,t)=g_{2}(x,t),\,\, x\in[0,1],\,\, t\in(0,0.5],\\
\end{split}
\end{equation}

 and the corresponding forcing terms $g_{1}(x,t)$ and $g_{2}(x,t)$  are of the form
\begin{equation}\label{91}
\begin{split}
&g_{1}(x,t)=e^{-it}\bigg(iu_{1}(x,0)-\lambda_{1}(-\Delta)^{\frac{\alpha}{2}}u_{1}(x,0)+2u_{1}(x,0)+u_{2}(x,0)\bigg),\\
&g_{2}(x,t)=e^{-it}\bigg(iu_{2}(x,0)-\lambda_{1}(-\Delta)^{\frac{\alpha}{2}}u_{2}(x,0)+2u_{2}(x,0)-u_{1}(x,0)\bigg).
\end{split}
\end{equation}
\end{exmp}
The exact solutions $u_{1}(x,t)=e^{-it}x^{7}$ and $u_{2}(x,t)=e^{-it}x^{7}$
with $\nu=1.1,\,\lambda_{1}=\frac{\Gamma(8-\nu)}{\Gamma(8)}$, $\lambda_{2}=\frac{\Gamma(8-\nu)}{\Gamma(8)}$. The errors and order of convergence are listed in Tables \ref{Tab:d1} and \ref{Tab:d2}, confirming optimal $O(h^{N+1})$ order of convergence across.

\begin{table}[!htb]
    \centering
\begin{center}
 \begin{tabular}{|c|| c c||c|| c c ||c||c c |}
  \hline
 \hline

 \hline
  N&\multicolumn{8}{|c|}{   N=1\quad \quad \quad\quad\quad \quad\qquad \qquad N=2\quad\quad\quad \qquad\quad \quad  \qquad \quad N=3}\\
\hline\hline
  K & $L^{2}$-Error & order & K &$L^{2}$-Error & order&K&$L^{2}$-Error& order \\ [0.5ex]

 \hline
92 & 2.27 e-02 & -&       60&1.93e-04&-&         50&4.1e-06&-\\

100 & 1.99e-02 &1.54&   90&5.60e-05&3.01&    70&1.23e-06&3.58\\
 130&1.07e-02&2.37&     110&3.0e-05&3.12&    100&2.98e-07&3.96\\

 \hline
 \hline

\end{tabular}
\end{center}
\caption{ $L^{2}$-Error and order of convergence for $u_{1}$ with $K$ elements and polynomial order $N$.}\label{Tab:d1}
\end{table}
\begin{table}[!htb]
    \centering
\begin{center}
 \begin{tabular}{|c|| c c||c|| c c ||c||c c |}
  \hline
 \hline

 \hline
  N&\multicolumn{8}{|c|}{   N=1\quad \quad \quad\quad\quad \quad\qquad \qquad N=2\quad\quad\quad \qquad\quad \quad  \qquad \quad N=3}\\
\hline\hline
  K & $L^{2}$-Error & order & K &$L^{2}$-Error & order&K&$L^{2}$-Error& order \\ [0.5ex]

 \hline
92 & 2.25e-02 & -&       60&1.7481e-04&-&         50&3.87e-06&-  \\

100 & 1.92e-02 &1.9&   90&5.03e-05&3.07&    70&8.91e-07&4.37\\
 130& 1.12e-02&2.04&     110&2.67e-05&3.16&    100&2.4e-07&3.68\\

 \hline
 \hline

\end{tabular}
\end{center}
\caption{ $L^{2}$-Error and order of convergence for $u_{2}$ with $K$ elements and polynomial order $N$.}\label{Tab:d2}
\end{table}
\begin{exmp}\label{ex6} We consider the nonlinear  coupled fractional Schr\"{o}dinger equations
\begin{equation}\label{sch1}
\begin{split}
&i\frac{\partial u_{1}(x,t)}{\partial t}- \lambda_{1}(-\Delta)^{\frac{\alpha}{2}}u_{1}(x,t)+u_{2}(x,t)+u_{1}(x,t)+(|u_{1}(x,t)|^{2}+|u_{2}(x,t)|^{2})
u_{1}(x,t)=g_{1}(x,t),\,\, x\in[0,1],\,\,t\in(0,0.5],\\\\
&i\frac{\partial u_{2}(x,t)}{\partial t}- \lambda_{2}(-\Delta)^{\frac{\alpha}{2}}u_{2}(x,t)+u_{2}(x,t)+u_{1}(x,t)+(|u_{1}(x,t)|^{2}
+|u_{2}(x,t)|^{2})u_{2}(x,t)=g_{2}(x,t),\, x\in[0,1],\, t\in(0,0.5],\\
\end{split}
\end{equation}

 and the corresponding forcing terms $g_{1}(x,t)$ and $g_{2}(x,t)$  are of the form
\begin{equation}\label{91}
\begin{split}
&g_{1}(x,t)=e^{-it}\bigg(iu_{1}(x,0)-\lambda_{1}(-\Delta)^{\frac{\alpha}{2}}u_{1}(x,0)+u_{2}(x,0)+u_{1}(x,0)
+(|u_{1}(x,0)|^{2}+|u_{1}(x,0)|^{2})u_{1}(x,0)\bigg),\\
&g_{2}(x,t)=e^{-it}\bigg(iu_{2}(x,0)-\lambda_{1}(-\Delta)^{\frac{\alpha}{2}}u_{2}(x,0)+u_{2}(x,0)+u_{1}(x,0)
+(|u_{1}(x,0)|^{2}+|u_{1}(x,0)|^{2})u_{2}(x,0)\bigg),
\end{split}
\end{equation}
\end{exmp}
to obtain an exact solutions $u_{1}(x,t)=e^{-it}x^{7}$ and $u_{2}(x,t)=e^{-it}x^{7}$
with $\nu=1.2,\,\lambda_{1}=\frac{\Gamma(8-\nu)}{2\Gamma(8)}$, $\lambda_{2}=\frac{\Gamma(8-\nu)}{2\Gamma(8)}$. The errors
and order of convergence are listed in Tables \ref{Tab:e1} and \ref{Tab:e2}, confirming optimal $O(h^{N+1})$ order of
convergence across.

\begin{table}[!htb]
    \centering
\begin{center}
 \begin{tabular}{|c|| c c||c|| c c ||c||c c |}
  \hline
 \hline

 \hline
  N&\multicolumn{8}{|c|}{   N=1\quad \quad \quad\quad\quad \quad\qquad \qquad N=2\quad\quad\quad \qquad\quad \quad  \qquad \quad N=3}\\
\hline\hline
  K & $L^{2}$-Error & order & K &$L^{2}$-Error & order&K&$L^{2}$-Error& order \\ [0.5ex]

 \hline
96 & 1.90 e-02& -&       30&4.7e-04&-&         40&8.68e-06&-  \\

120 & 1.27e-02 &2.35&   60&1.47e-04&2.86 &    60&1.79e-06&3.89 \\
 135&9.6e-03&1.92&     130&1.22e-05&3.22&    80&6.03e-07&3.78\\

 \hline
 \hline

\end{tabular}
\end{center}
\caption{ $L^{2}$-Error and order of convergence for $u_{1}$ with $K$ elements and polynomial order $N$.}\label{Tab:e1}
\end{table}

\begin{table}[!htb]
    \centering
\begin{center}
 \begin{tabular}{|c|| c c||c|| c c ||c||c c |}
  \hline
 \hline

 \hline
  N&\multicolumn{8}{|c|}{   N=1\quad \quad \quad\quad\quad \quad\qquad \qquad N=2\quad\quad\quad \qquad\quad \quad  \qquad \quad N=3}\\
\hline\hline
  K & $L^{2}$-Error & order & K &$L^{2}$-Error & order&K&$L^{2}$-Error& order \\ [0.5ex]

 \hline
96 & 1.89e-02 & -&       40&4.18e-04&-&         40&7.71e-06&-  \\

120 & 1.34e-02 &1.55&   60&1.26e-04&2.95 &    60&1.47e-06&4.08\\
 135&1.03e-02&2.22&     130&1.21e-05&3.04&    90&5.1e-07&3.7\\

 \hline
 \hline

\end{tabular}
\end{center}
\caption{ $L^{2}$-Error and order of convergence for $u_{2}$ with $K$ elements and polynomial order $N$.}\label{Tab:e2}
\end{table}

\begin{figure}
\begin{center}
  \vspace{4mm}
\begin{tabular}{cc}
\hspace{-1.0cm}
\begin{overpic}[scale=0.35]{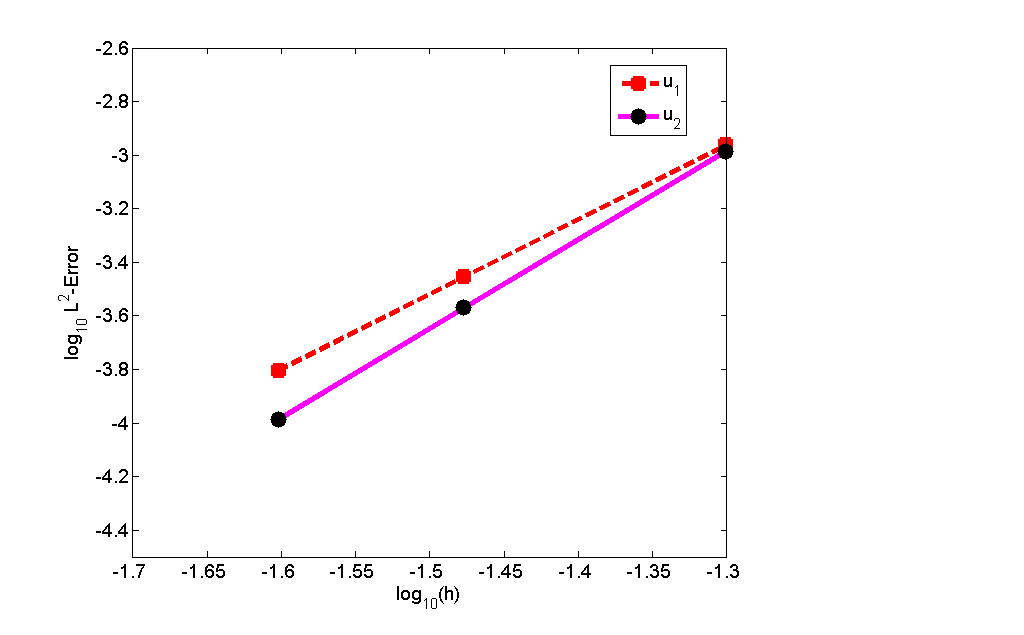}
\end{overpic}
 &\hspace{-1.8cm}
\begin{overpic}[scale=0.35]{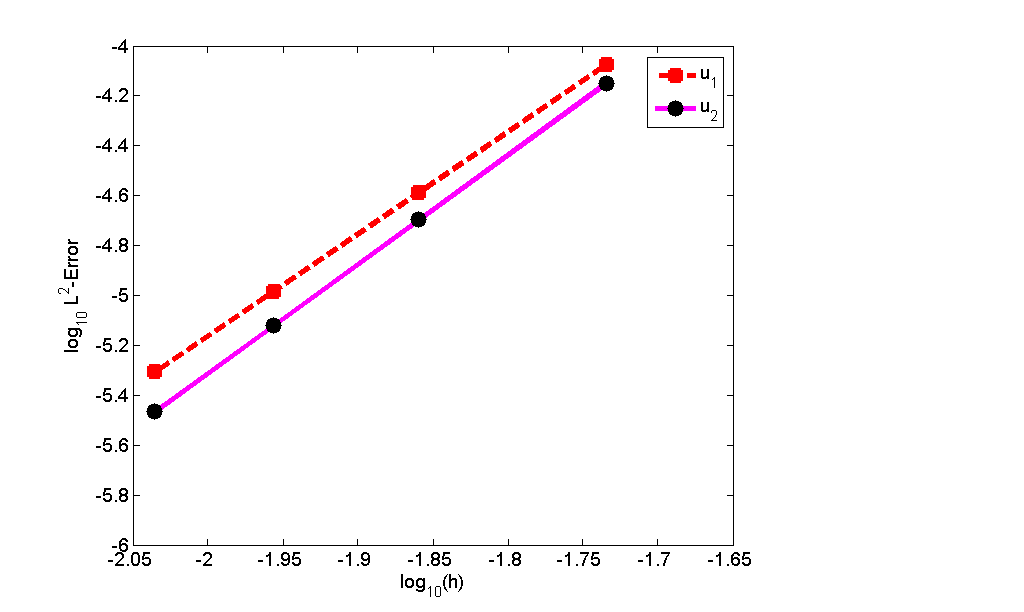}
\end{overpic}
\end{tabular}
\end{center}
\caption[]{\small{ The convergence rate  of \eqref{sch1mn2} for  $N = 2$ (left), $N = 3$ (right).}} \label{fig:h}
\end{figure}

\begin{exmp}\label{ex12vn} Consider  the following nonlinear  coupled fractional Schr\"{o}dinger equations
\begin{equation}\label{sch1mn2}
\begin{split}
&i\frac{\partial u_{1}(x,t)}{\partial t}- \lambda_{1}(-\Delta)^{\frac{\alpha}{2}}u_{1}(x,t)+u_{2}(x,t)+u_{1}(x,t)+(|u_{1}(x,t)|^{2}
+|u_{2}(x,t)|^{2})u_{1}(x,t)=g_{1}(x,t),\, x\in[-1,1],\, t\in(0,0.5],\\\\
&i\frac{\partial u_{2}(x,t)}{\partial t}- \lambda_{2}(-\Delta)^{\frac{\alpha}{2}}u_{2}(x,t)+u_{2}(x,t)-u_{1}(x,t)+(|u_{1}(x,t)|^{2}+|u_{2}(x,t)|^{2})
u_{2}(x,t)=g_{2}(x,t),\, x\in[-1,1],\, t\in(0,0.5],\\
\end{split}
\end{equation}
and the corresponding forcing terms $g_{1}(x,t)$ and $g_{2}(x,t)$  are of the form
\begin{equation}\label{91}
\begin{split}
&g_{1}(x,t)=e^{-it}\bigg(iu_{1}(x,0)-\lambda_{1}(-\Delta)^{\frac{\alpha}{2}}u_{1}(x,0)+u_{2}(x,0)+u_{1}(x,0)
+(|u_{1}(x,0)|^{2}+|u_{1}(x,0)|^{2})u_{1}(x,0)\bigg),\\
&g_{2}(x,t)=e^{-it}\bigg(iu_{2}(x,0)-\lambda_{1}(-\Delta)^{\frac{\alpha}{2}}u_{2}(x,0)+u_{2}(x,0)-u_{1}(x,0)
+(|u_{1}(x,0)|^{2}+|u_{1}(x,0)|^{2})u_{2}(x,0)\bigg),
\end{split}
\end{equation}
\end{exmp}
 The exact solutions $u_{1}(x,t)=e^{-it}(x^{2}-1)^{6}$ and $u_{2}(x,t)=e^{-it}(x^{2}-1)^{6}$
with $\nu=1.3,\,\lambda_{1}=\frac{\Gamma(13-\nu)}{2\Gamma(13)}$, $\lambda_{2}=\frac{\Gamma(13-\nu)}{2\Gamma(13)}$. We consider cases with $N = 2, 3$ and $\log_{10}(h)$. The numerical orders of convergence
are shown in Figure \ref{fig:h}, showing an
$O(h^{N+1})$ convergence rate for all orders.

\begin{figure}
\begin{center}
  \vspace{6mm}
\begin{tabular}{cc}
\hspace{-6cm}
\begin{overpic}[scale=0.4]{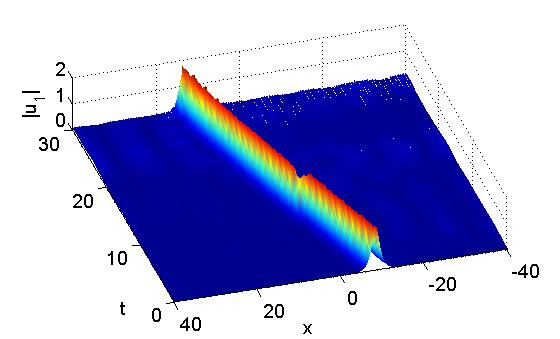}
\end{overpic}
 &\hspace{-4.2cm}
\begin{overpic}[scale=0.4]{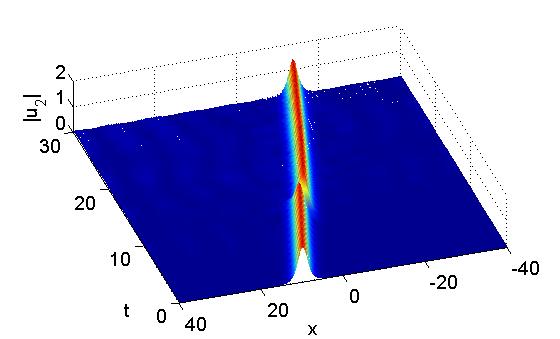}
\end{overpic}
\\
\\
\\
\hspace{1cm}
\begin{overpic}[scale=0.4]{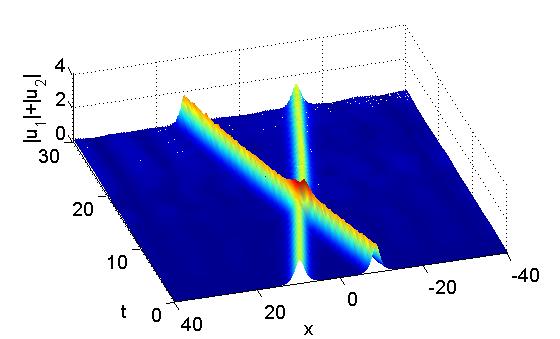}
\end{overpic}

\end{tabular}
\caption[]{{  Numerical solutions for Example \ref{ex12vng123}  with $\beta= 1$ and $\alpha = 2$.}}\label{fig:elastic1}
\end{center}
\end{figure}

\begin{figure}
\begin{center}
  \vspace{6mm}
\begin{tabular}{cc}
\hspace{-6cm}
\begin{overpic}[scale=0.4]{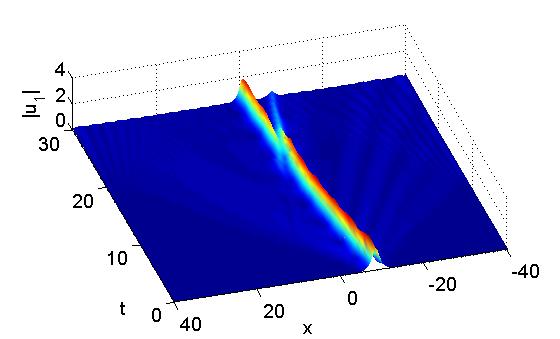}
\end{overpic}
 &\hspace{-4.0cm}
\begin{overpic}[scale=0.4]{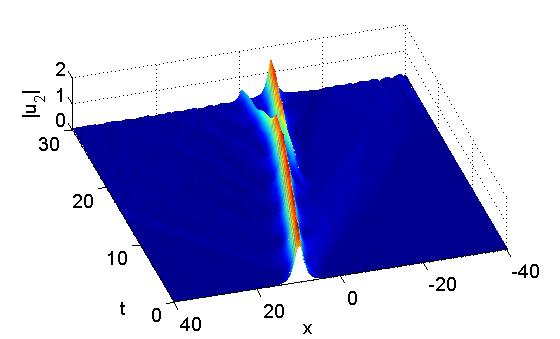}
\end{overpic}
\\
\\
\\
\hspace{1cm}
\begin{overpic}[scale=0.4]{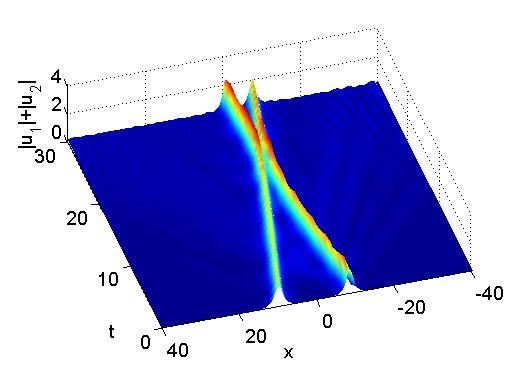}
\end{overpic}

\end{tabular}
\caption[]{{  Numerical solutions for Example \ref{ex12vng123}  with $\beta= 1$ and $\alpha = 1.6$.}}\label{fig:elastic2}
\end{center}
\end{figure}

\begin{figure}
\begin{center}
  \vspace{6mm}
\begin{tabular}{cc}
\hspace{-6cm}
\begin{overpic}[scale=0.4]{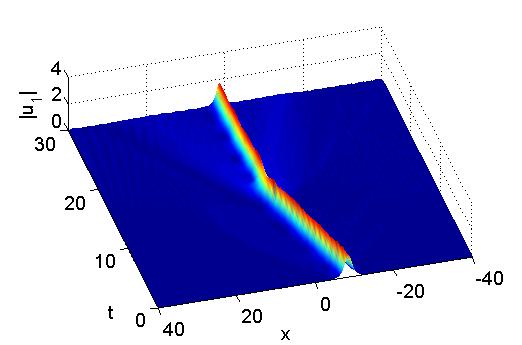}
\end{overpic}
 &\hspace{-4.0cm}
\begin{overpic}[scale=0.4]{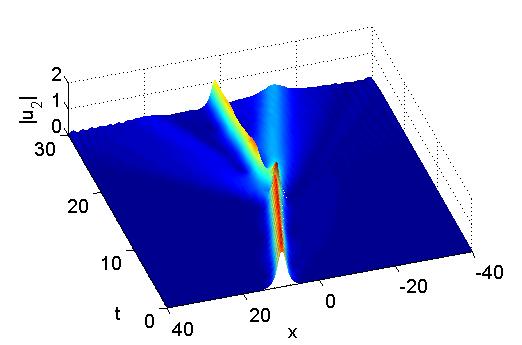}
\end{overpic}
\\
\\
\\
\hspace{1cm}
\begin{overpic}[scale=0.4]{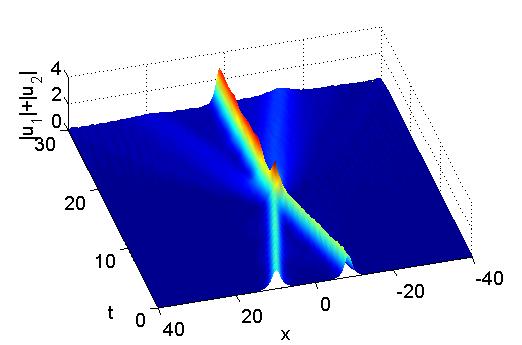}
\end{overpic}

\end{tabular}
\caption[]{{  Numerical solutions for Example \ref{ex12vng123}  with $\beta= 0.3$ and $\alpha = 1.8$.}}\label{fig:elastic2nb}
\end{center}
\end{figure}

\begin{exmp}\label{ex12vng123} We consider the following weakly coupled problem
\begin{equation}\label{sch1mn2}
\begin{split}
&i\frac{\partial u_{1}}{\partial t}- (-\Delta)^{\frac{\alpha}{2}}u_{1}+(|u_{1}|^{2}+\beta|u_{2}|^{2})u_{1}=0,\\
&i\frac{\partial u_{2}}{\partial t}- (-\Delta)^{\frac{\alpha}{2}}u_{2}+(\beta|u_{1}|^{2}+|u_{2}|^{2})u_{2}=0,\\
\end{split}
\end{equation}
 subject to the initial conditions

\begin{equation}\label{91}
\begin{split}
u_{1}(x,0)=\sqrt{2}r_{1}sech(r_{1}x+D)e^{iV_{0}x},\\
u_{2}(x,0)=\sqrt{2}r_{2}sech(r_{2}x+D)e^{iV_{0}x},\\
\end{split}
\end{equation}
when $\beta= 1$ and $\alpha = 2$, the problem collapses to the Manakov equation, and
the solitary waves collide elastically see Figure \ref{fig:elastic1}. The exact solutions are given by
\begin{equation}\label{91}
\begin{split}
u_{1}(x,t)=\sqrt{2}r_{1}sech(r_{1}x-2r_{1}V_{0}t+D)e^{i(V_{0}x+(r_{1}^{2}-V_{0}^{2})t)},\\
u_{2}(x,t)=\sqrt{2}r_{2}sech(r_{2}x-2r_{2}V_{0}t-D)e^{i(-V_{0}x+(r_{2}^{2}-V_{0}^{2})t)},\\
\end{split}
\end{equation}
\end{exmp}
where $r_{1} = 1$, $r_{2} = 1$, $V_{0}=0.4$, $D=10$ and $x\in[-40,40]$. The Figures \ref{fig:elastic2} and \ref{fig:elastic2nb} present the numerical solutions for different values of order $\alpha$ and $\beta$.  From these figures  it is obvious that  the collision of solitons are inelastic. In particular, the colliding particles stick together after interaction when $\alpha =1.8$,
which means that there may occur a completely inelastic collision see Figure \ref{fig:elastic2nb}.

\begin{figure}
\begin{center}
  \vspace{6mm}
\begin{tabular}{cc}
\hspace{-6cm}
\begin{overpic}[scale=0.4]{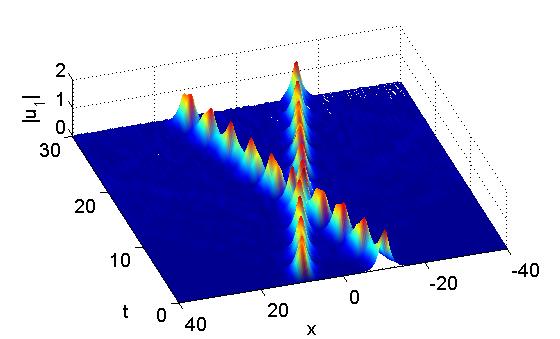}
\end{overpic}
 &\hspace{-4.0cm}
\begin{overpic}[scale=0.4]{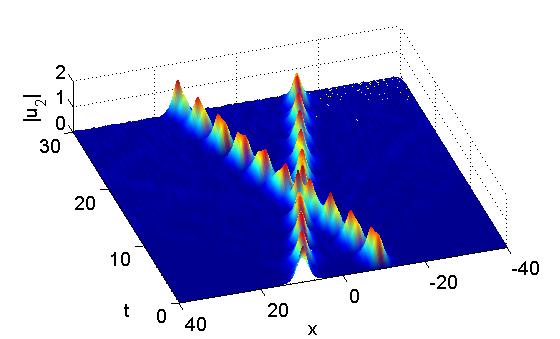}
\end{overpic}
\\
\\
\\
\hspace{1cm}
\begin{overpic}[scale=0.4]{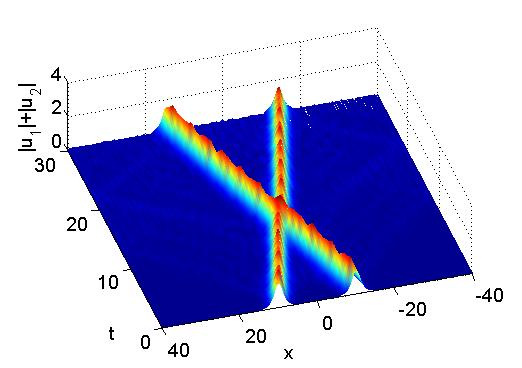}
\end{overpic}

\end{tabular}
\caption[]{{  Numerical solutions for Example \ref{ex12vng} with $\varpi_{1}= 1$ , $\alpha = 2$.}}\label{fig:elastic3}
\end{center}
\end{figure}

\begin{figure}
\begin{center}
  \vspace{6mm}
\begin{tabular}{cc}
\hspace{-6cm}
\begin{overpic}[scale=0.4]{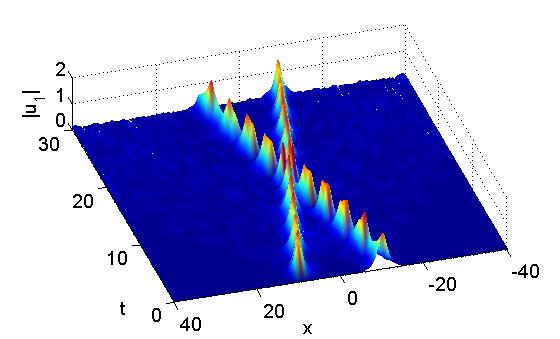}
\end{overpic}
 &\hspace{-4.0cm}
\begin{overpic}[scale=0.4]{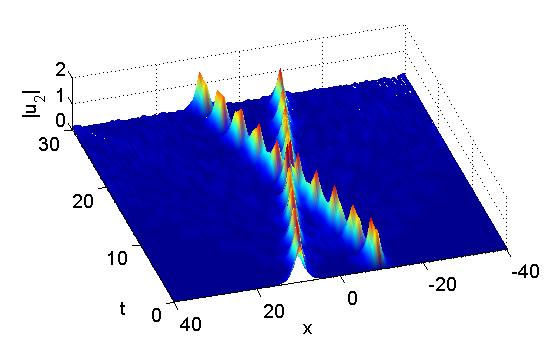}
\end{overpic}
\\
\\
\\
\hspace{1cm}
\begin{overpic}[scale=0.4]{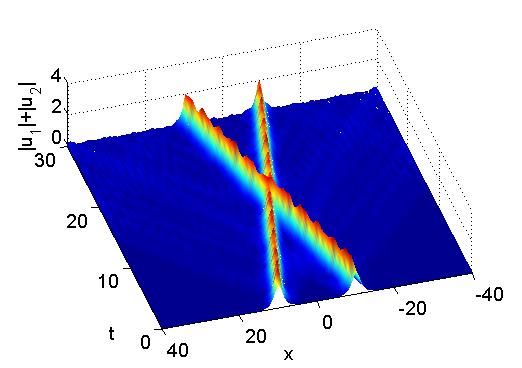}
\end{overpic}

\end{tabular}
\caption[]{{  Numerical solutions for Example \ref{ex12vng} with $\varpi_{1}= 1$ , $\alpha = 1.6$. }}\label{fig:elastic4}
\end{center}
\end{figure}

\begin{figure}
\begin{center}
  \vspace{6mm}
\begin{tabular}{cc}
\hspace{-6cm}
\begin{overpic}[scale=0.4]{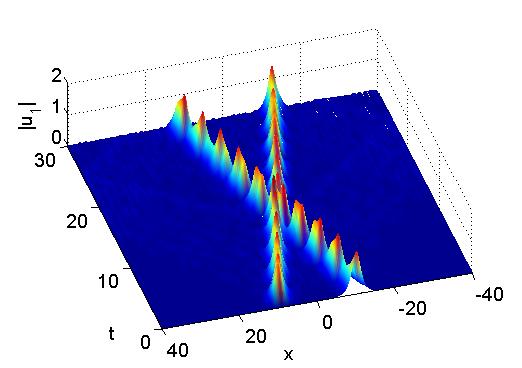}
\end{overpic}
 &\hspace{-4.0cm}
\begin{overpic}[scale=0.4]{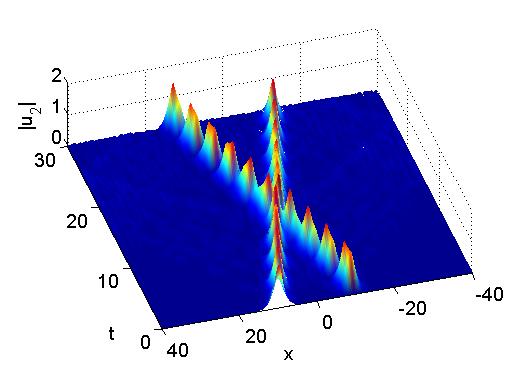}
\end{overpic}
\\
\\
\\
\hspace{1cm}
\begin{overpic}[scale=0.4]{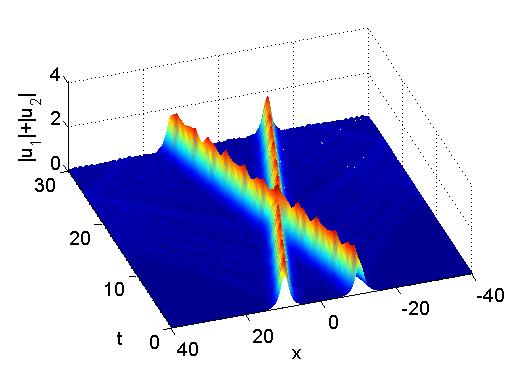}
\end{overpic}

\end{tabular}
\caption[]{{  Numerical solutions for Example \ref{ex12vng}  with $\varpi_{1}= 1$ and $\alpha = 1.8$.}}\label{fig:elastic223}
\end{center}
\end{figure}

\begin{figure}
\begin{center}
  \vspace{6mm}
\begin{tabular}{cc}
\hspace{-6cm}
\begin{overpic}[scale=0.4]{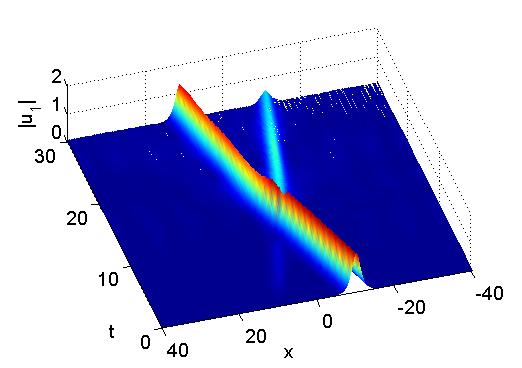}
\end{overpic}
 &\hspace{-4.0cm}
\begin{overpic}[scale=0.4]{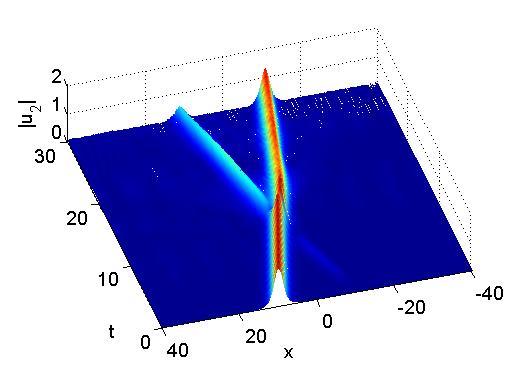}
\end{overpic}
\\
\\
\\
\hspace{1cm}
\begin{overpic}[scale=0.4]{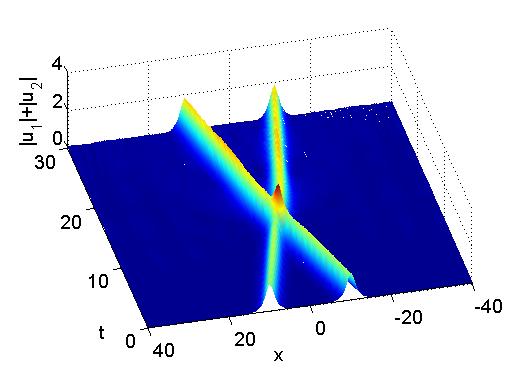}
\end{overpic}

\end{tabular}
\caption[]{{  Numerical solutions for Example \ref{ex12vng}  with $\varpi_{1}= 0.0175$ and $\alpha = 2$.}}\label{fig:inelastic6}
\end{center}
\end{figure}

\begin{figure}
\begin{center}
  \vspace{6mm}
\begin{tabular}{cc}
\hspace{-1cm}
\begin{overpic}[scale=0.30]{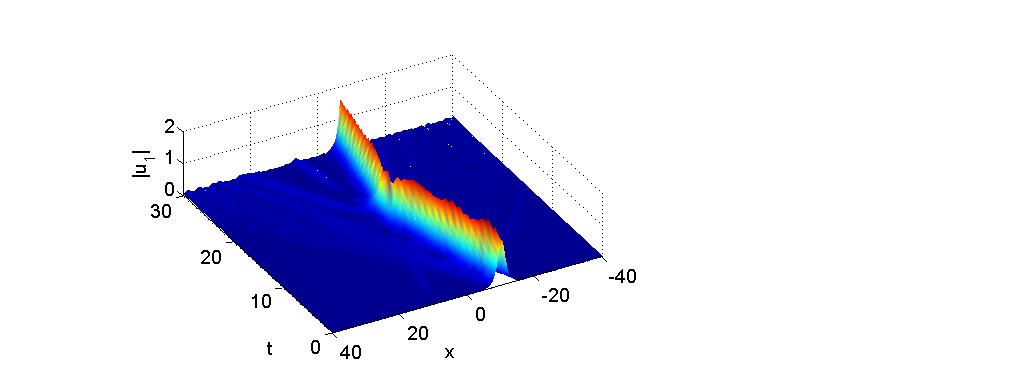}
\end{overpic}
 &\hspace{-4.6cm}
\begin{overpic}[scale=0.30]{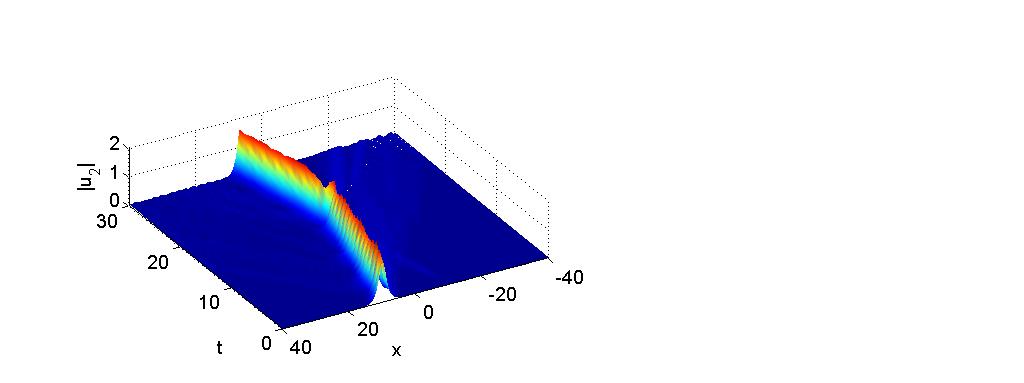}
\end{overpic}
\\
\\
\\
\hspace{2cm}
\begin{overpic}[scale=0.4]{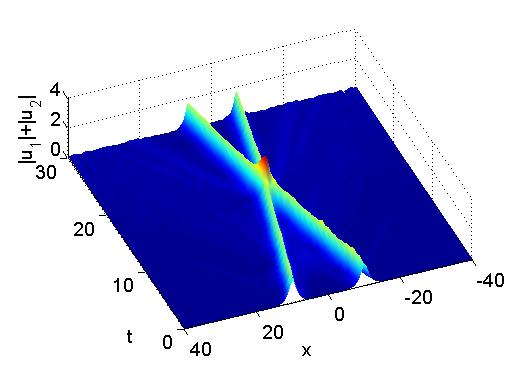}
\end{overpic}

\end{tabular}
\caption[]{{  Numerical solutions for Example \ref{ex12vng}  with $\varpi_{1}= 0.0175$ and $\alpha = 1.6$.}}\label{fig:inelastic6h}
\end{center}
\end{figure}

\begin{figure}
\begin{center}
  \vspace{6mm}
\begin{tabular}{cc}
\hspace{-6cm}
\begin{overpic}[scale=0.4]{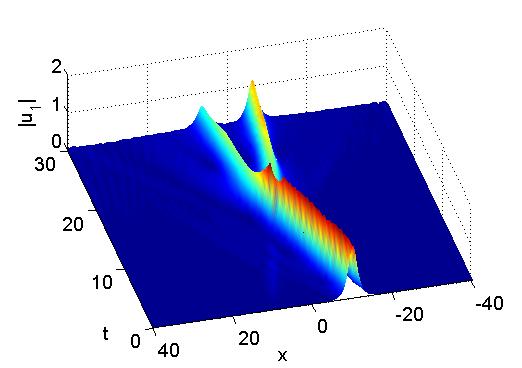}
\end{overpic}
 &\hspace{-4.0cm}
\begin{overpic}[scale=0.4]{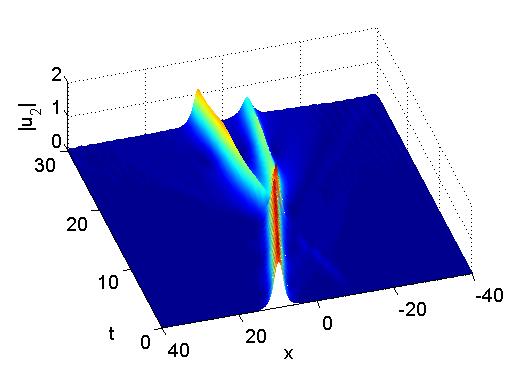}
\end{overpic}
\\
\\
\\
\hspace{1cm}
\begin{overpic}[scale=0.4]{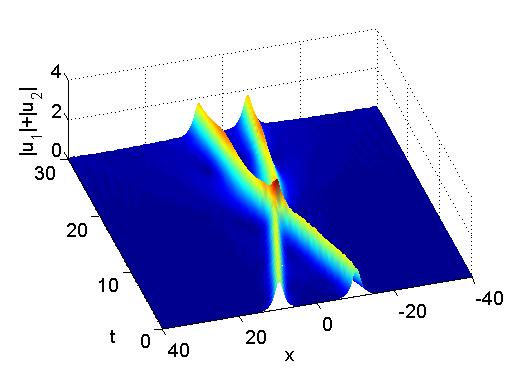}
\end{overpic}

\end{tabular}
\caption[]{{  Numerical solutions for Example \ref{ex12vng} with $\varpi_{1}= 0.0175$ and $\alpha = 1.8$.}}\label{fig:inelastic8h}
\end{center}
\end{figure}

\begin{exmp}\label{ex12vng} Finally,  we consider the strongly coupled system as follows

\begin{equation}\label{sch1mn2}
\begin{split}
&i\frac{\partial u_{1}}{\partial t}- (-\Delta)^{\frac{\alpha}{2}}u_{1}+(|u_{1}|^{2}+|u_{2}|^{2})u_{1}+u_{1}+\varpi_{1}u_{2}=0,\\
&i\frac{\partial u_{2}}{\partial t}- (-\Delta)^{\frac{\alpha}{2}}u_{2}+(|u_{1}|^{2}+|u_{2}|^{2})u_{2}+\varpi_{1}u_{1}+u_{2}=0,\\
\end{split}
\end{equation}
 subject to the initial conditions

\begin{equation}\label{91}
\begin{split}
u_{1}(x,0)=\sqrt{2}r_{1}sech(r_{1}x+D)e^{iV_{0}x},\\
u_{2}(x,0)=\sqrt{2}r_{2}sech(r_{2}x+D)e^{iV_{0}x},\\
\end{split}
\end{equation}

\end{exmp}
where $r_{1} =r_{2} = 1$, $V_{0}=0.4$, $D=10$ and $x\in[-40,40]$. \\
Elastic collisions: The collision of the solitary waves is elastic \cite{aydin2011lobatto} when $\varpi_{1}= 1$, $\alpha = 2$ see Figure \ref{fig:elastic3}. We observe that the two waves emerge without any changes in their shapes and velocities after collision. Taking $\varpi_{1}= 1$,  we compute the numerical solutions for different values of $\alpha$, which are depicted in Figures \ref{fig:elastic4} and \ref{fig:elastic223}.  From these figures,  for any $1 < \alpha \leq 2$, the collision is always elastic. When $\alpha$ tends to $2$, the shape of the solitons will change more slightly and the waveforms become closer to the classical case with $\alpha = 2$.\\
Inelastic collision: The collision is inelastic \cite{aydin2011lobatto} when $\varpi_{1}= 0.0175$ and $\alpha = 2$ see Figure \ref{fig:inelastic6}. It is clear that the shapes and directions of two waves have changed after interaction. The observation is in accordance with the known result.\\
 The Figures \ref{fig:inelastic6h} and \ref{fig:inelastic8h} present the numerical solutions for different values of order $\alpha$ for fixed $\varpi_{1}= 0.0175$.  From these figures it is obvious that the collision is always inelastic.


\newpage

\section{Conclusions}\label{sc6}
In this work, we developed and analyzed a nodal discontinuous Galerkin  method for solving the   nonlinear
 fractional Schr\"{o}dinger equation and the strongly coupled nonlinear fractional Schr\"{o}dinger equations, and have proven the stability of these methods. They are  discretized  using  high-order nodal  basis set of orthonormal Lagrange-Legendre polynomials as a more suitable and computationally stable approach. Numerical experiments confirm that the optimal order of convergence is recovered. As a last two examples, the  weakly coupled nonlinear fractional Schr\"{o}dinger equations  with initial conditions are solved for different values of $\alpha$ and results show that  the collision of solitons are inelastic when $\alpha\neq2$ and  the results of the strongly
  nonlinear fractional Schr\"{o}dinger equations  are the shape of the soliton will change  slightly as $\alpha$ increase, with the classical case $\varpi_{1}= 1$ and $\alpha =2$ as the limit. When $\varpi_{1}= 1$ and $ \alpha\neq2$,  the collision is always elastic and the collision is inelastic when $\varpi_{1}= 0.0175$ and $1<\alpha \leq 2$.
 \\
%


\end{document}